\documentclass[a4,12pt,multicol,makeidx]{article}
\def\origin{%
\hbox{}\vskip-\baselineskip\vskip-\topskip%
  \vbox to 0pt{\vskip-1in%
    \hbox to 0pt{\hskip-1in%
      \hbox to 0pt{\vrule width 1cm height .4pt depth 0mm\hss}%
      \vbox to 0pt{\hrule width .4pt height 0pt depth 1cm\vss}%
    \hss}%
  \vss}
  \vskip-\baselineskip
  \vbox to 0pt{\vskip-1in\vskip3cm%
    \hbox to 0pt{\hskip-1in\hskip3cm%
      \hbox to 0pt{\hss\vrule width 2cm height .4pt depth 0mm\hss}%
      \vbox to 0pt{\vss\hrule width .4pt height 1cm depth 1cm\vss}%
    \hss}%
  \vss}%
  \vskip5mm\hskip10mm (3cm,3cm)
}%
 \def\a{\alpha} \def\b{\beta}  \def\l{\lambda}   \def\p{\partial}  \def\e{\varepsilon} 
\def\n{\nabla}    
\def\leq{\underline{<}} \def\la{\langle} \def\ra{\rangle}
\def\hx{\hat{x}} \def\hy{\hat{y}}  
\def\ox{\overline{x}} \def\oy{\overline{y}}   
\def\ov{\overline{v}} \def\uv{\underline{v}} \def\ou{\overline{u}} 
\def\oI{\overline{I}}

\newenvironment{theorem}{%
\par \bigskip \it}{%
\bigskip \par}

\newenvironment{definition}{%
\par \bigskip \it}{%
\bigskip \par}

\makeindex
\title{Homogenizations of integro-differential equations with L{\'e}vy operators with asymmetric and degenerate  densities.\\
}
\author{Mariko Arisawa\\ DAMTP, Centre for Mathematical Sciences\\
University of Cambridge\\
Wilberforce road\\
Cambridge, CB3 0WA, England\\
E-mail: M.Arisawa@damtp.cam.ac.uk
}
\date{}
\begin{document}
\maketitle
\bigskip


{\bf Synopsis.} We consider periodic homogenization problems for the L{\'e}vy operators with asymmetric L{\'e}vy densities. 
The formal asymptotic expansion used for the $\a$-stable (symmetric) L{\'e}vy operators ($\a\in (0,2)$) is not applicable directly 
 to such asymmetric cases. We rescale the asymmetric densities, extract the most singular part of the measures, which average 
out the microscopic dependences in the homogenization procedures. We give two conditions (A) and (B), which characterize such a 
class of asymmetric densities, that the above "rescaled" homogenization is available. \\ 

\section{Introduction.} 

$\quad$In this paper, we are interested in the following homogenization problems concerning with the L{\'e}vy operator : 
\begin{equation}\label{homo}
	u_{\e}(x)
	-a(\frac{x}{\e})\int_{{\bf R^M}}  [u_{\e}(x+\beta(z))-u_{\e}(x)- {\bf 1}_{|z|\leq 1}\la \n u_{\e}(x),\beta(z) \ra] dq(z) 
\end{equation} 
$$\qquad\qquad\qquad\qquad\qquad\qquad\qquad\qquad\qquad\qquad
	-f(\frac{x}{\e})=0\quad \hbox{in}\quad \Omega, 
$$
\begin{equation}\label{bc}
	u=\phi(x)\qquad\qquad\qquad\qquad\qquad\qquad\qquad\qquad\qquad\qquad\quad\quad \hbox{in}\quad \Omega^c, 
\end{equation}
and 
\begin{equation}\label{homo2}
	u_{\e}(x)
	-a(\frac{x}{\e})\int_{{\bf R^M}}  [u_{\e}(x+\beta(z))-u_{\e}(x)] dq(z) -f(\frac{x}{\e})=0\quad \hbox{in}\quad \Omega, 
\end{equation} 
with (\ref{bc}). Here, 
$\Omega$ is an open bounded domain in ${\bf R^N}$, $M\leq N$, 
$\beta$ is a positively homogenious, continuous function from ${\bf R^M}$ to ${\bf R^N}$ such that
\begin{equation}\label{beta}
	\beta(cz)=c\beta(z)\quad \forall c>0;      \quad |\beta(z)| \leq B_1|z|\qquad \forall z\in {\bf R^M},
\end{equation}
where $B_1>0$ is a constant, 
$dq(z)=q(z)dz$ is a positive Radon measure on ${\bf R^M}$ which satisfies
\begin{equation}\label{radon}
	\int_{|z|< 1} |z|^{\gamma} dq(z)+\int_{|z| \geq 1} |z|^{\gamma -1} dq(z) <\infty,
\end{equation}
with $\gamma=2$ in the case of (\ref{homo}), and with $\gamma=1$ in the case of (\ref{homo2}), 
 $a$, $f$ are real valued continuous functions defined in ${\bf R^N}$, periodic in ${\bf T^N}$$=[0,1]^N$, such that there exist constants $\theta_1$, $\theta_2\in (0,1]$, $L>0$, $a_0>0$ with which the following hold:
\begin{equation}\label{holder1}
	 a(\cdot)\geq \exists a_0>0; \quad |a(y)-a(y')|\leq L|y-y'|^{\theta_1} \quad y,y'\in {\bf R^N}, 
\end{equation}
\begin{equation}\label{holder2}
	|f(y)-f(y')|\leq L|y-y'|^{\theta_2} \quad y,y'\in {\bf R^N},
\end{equation}
and $\phi$ is a real valued bounded continuous function defined in $\Omega^c$. \\

For any $\e>0$, there exists a unique solution  $u_{\e}$ of (\ref{homo})-(\ref{bc}) and of (\ref{homo2})-(\ref{bc}) respectively  in the framework of the viscosity solution (see \S 6 for the definition, and also  M. Arisawa \cite{arnewdef}, \cite{arcorrig}, \cite{arcomp}, G. Barles and C. Imbert \cite{bi} for the existence and the uniqueness results, and 
M.G. Crandall, H. Ishii, and P.-L. Lions \cite{users} for the general theory of the viscosity solution).  As $\e$ goes to zero the sequence of functions $\{u_{\e}\}$ converges locally uniformly to a limit $\ou$, and we are interested in 
 finding an effective nonlocal equation which characterizes $\ou$. 

Such a homogenization problem was solved in the case that the L{\'e}vy measure is $\a$-stable (see M. Arisawa \cite{arperiodic}, \cite{arquasi}):  
$$
	dq(z)=\frac{1}{|z|^{N+\a}} \qquad z\in {\bf R^M}, \quad \a\in (0,2)\quad \hbox{a fixed number},
$$
 by utilizing the formal asymptotic expansion : 
\begin{equation}\label{f1}
	u_{\e}(x)=\overline{u}(x)+\e^{\a}v(\frac{x}{\e}) + o(\e^{\a}) \qquad x\in {\bf R^N},
\end{equation}
where $\overline{u}=\lim_{\e\to 0} u_{\e}$, $v$  a periodic function defined in ${\bf R^N}$, called corrector. The above expansion leads to the so-called  the ergodic cell problem, which gives the effective equation for $\overline{u}$. We refer the readers to 
A. Bensoussan, J.L. Lions, and G. Papanicolaou \cite{blp} for the detailed discussion on this method.  In the framework of the viscosity solution, the formal argument can be  justified rigorously, by the perturbed test function method established by L.C. Evans in \cite{ev1} and \cite{ev2} (see also P.-L. Lions, G. Papanicolaou, and S.R.S. Varadhan \cite{lpv}). However, as we shall see below in  Examples 1-4, 
the above formal expansion cannot be employed directly if the measure $dq(z)$ is asymmetric. Here, we assume that the L{\'e}vy measure 
satisfies the following condition (A).\\

{\bf (A)}:  Let $S=\hbox{supp}(dq(z))\subset {\bf R^M}$. There exists  a constant  $\a\in (0,2)$ such that 
\begin{equation}\label{bound}
	\e^{M+\a}q(\e z) \leq  C_1 |z|^{-(M+\a)}  \quad \forall \e\in (0,1), \quad  \forall z\in {\bf R^M}, 
\end{equation}
 where $C_1>0$ is a constant independent on $\e$, 
 and a subset $S_0\subset S$ and a positive function $q_0(z)$ $(z\in {\bf R^M})$ such that 
\begin{equation}\label{lim}
	\lim_{\e \downarrow 0} \e^{M+\a} q(\e z) = q_0(z) \quad \forall z\in S_0; \quad =0 \quad \forall z\in {\bf R^M}\backslash S_0. 
\end{equation}
$$\quad$$
We define a new measure: 
\begin{equation}\label{q0}
	dq_{0}(z)=q_0(z)dz \quad \forall z\in S_0;\quad =0 dz\quad  \forall z\in  {\bf R^M}\backslash S_0. 
\end{equation}
The following property holds for this rescaled measure $dq_{0}(z)$. \\

{\bf Lemma 1.1.$\quad$}
\begin{theorem}  Assume that the Radon measure $dq(z)$ satisfies (\ref{radon}) and the condition (A). Then, $S_0$ is a positive cone, i.e. 
$$
	sS_0\subset S_0 \quad \forall s>0. 
$$
Moreover,  $s^{M+\a} q_0(sz)$$=q_0(z)$ ($\forall s>0$, $\forall z\in S_0$), and
\begin{equation}\label{power2}
	q_0(z) = |z|^{M+\a} \overline{q}_{0}(arg z)\quad \forall z\in {\bf R^M},
\end{equation}
where $\overline{q}_{0}(\theta)$ ($\theta\in [0,2\pi)$) is a bounded real valued function. 
\end{theorem} 

$Proof.$ Let $z\in S_0$. For any $s\in (0,1)$, from the condition (A), (\ref{lim}), 
$$
	\lim_{\e\to 0} \e^{M+\a} q(\e sz) = s^{-(M+\a)}\lim_{\e\to 0} (\e s)^{M+\a} q(\e sz) \qquad\qquad\qquad\qquad\qquad
$$
$$
	=s^{-(M+\a)}    \lim_{\e'\to 0} \e'^{M+\a} q(\e' z)  =
	s^{-(M+\a)}q_0(z)> 0. 
$$
Thus, $sz\in S_0$, and 
$$
	q_0(sz)= s^{-(M+\a)}q_0(z) \quad \forall s\in (0,1), \quad \forall z\in S_0. 
$$
Therefore, $q_0(z)=|z|^{M+\a}$$q_0(\frac{z}{|z|})$, and from the condition (A) (\ref{power2}) is proved. 
 \\

The following examples satisfy the condition (A).\\

{\bf Example 1.} Let $M=N$, $\beta(z)=z$, and for $\a\in (1,2)$ (resp. $(0,1)$), 
$$
	dq(z)=|z|^{-(M+\a)}dz \quad z\in {\bf R_{+}^{M}};\quad =0 \quad z\in ({\bf R_{+}^{M}})^c,
$$
where ${\bf R_{+}^{M}}=\{z=(z_1,...,z_M)| \quad z_i>0\quad 1\leq \forall i\leq M\}$. In this case, for $S=S_0={\bf R_{+}^{M}}$, we have
$$
	q(\e z)\e^{M+\a}=|z|^{-(M+\a)}=q_0(z)\quad \forall  z\in S_0; \quad =0 \quad \forall  z\in S_0^c,\quad \forall \e>0,
$$
and for 
$$
	dq_0(z)=|z|^{-(M+\a)}dz \quad z\in {\bf R_{+}^{M}};\quad =0 \quad z\in ({\bf R_{+}^{M}})^c, 
$$
the condition (A) is satisfied. Both $dq(z)$ and $dq_0(z)$ satisfy (\ref{radon}) with $\gamma=2$ (resp. $\gamma=1$).\\

{\bf Example 2.} Let $M=N=1$, $\beta(z)=z$, and for $1< \a_1<\a_2< 2$, 
$$
	dq(z)=|z|^{-(1+\a_1)}dz \quad z\leq -1 \quad \hbox{and}\quad z>0;\quad = |z|^{-(1+\a_2)} \quad -1<z<0.
$$
 In this case, for $\a=\a_2$, $S={\bf R}$, $S_0=\{z\in {\bf R}|\quad z<0\}$, we have  
$$
	\lim_{\e\to 0} q(\e z)\e^{1+\a}=|z|^{-(1+\a_2)}=q_0(z)\quad \forall z\in S_0; \quad =0\quad \forall z\in S_0^c,
$$
and for 
$$
	dq_0(z)=0 \quad z>0;\quad = |z|^{-(1+\a)} \quad z<0, 
$$
the condition (A) is satisfied. Both $dq(z)$ and $dq_0(z)$ satisfy (\ref{radon}) with $\gamma=2$.\\

{\bf Example 3.} Let $M=1$, $N=2$, $\beta(z)=(z,\gamma z)$, where $\gamma>0$ is an irrational number, and for $\a\in (1,2)$ (resp. ($0,1$)), 
$$
	dq(z)=|z|^{-(1+\a)}dz \quad z\in {\bf R}.
$$
In this case, for $S=S_0={\bf R}$, we have 
$$
	q(\e z)\e^{1+\a}=|z|^{-(1+\a)}=q_0(z) \quad  \forall z\in S=S_0,\quad \forall \e>0,
$$
and for 
$$
	dq_0(z)=|z|^{-(N+\a)}dz \quad z\in {\bf R}, 
$$
the condition (A) is satisfied. Both $dq(z)$ and $dq_0(z)$ satisfy (\ref{radon}) with $\gamma=2$ (resp. $\gamma=1$).\\

{\bf Example 4.} Let $M=N$, $\beta(z)=z$, and for $\gamma>0$, $1< \a < 2$, 
$$
	dq(z)=\exp(-\gamma|z|)|z|^{-(M+\a)}dz \quad z\in {\bf R^M}.
$$
In this case, for $g(s)=s^{\a}$, $S=S_0={\bf R^M}$, we have 
$$
	\lim_{\e\to 0} q(\e z)\e^{M+\a}=\lim_{\e\to 0} \exp(-\e \gamma|z|)   |z|^{-(M+\a)}  = |z|^{-(M+\a)}  \quad \forall z\in S=S_0,
$$
and for 
$$
	dq_0(z)=|z|^{-(M+\a)}dz \quad z\in {\bf R^M}, 
$$
the condition (A) is satisfied. Both $dq(z)$ and $dq_0(z)$ satisfy (\ref{radon}) with $\gamma=2$.\\

In Examples 1-3, the  L{\'e}vy measures are either asymmetric or degenerate (in the sense that $S$ or $S_0$ does not contain an open ball centered at the origin in ${\bf R^M}$). Example 3 corresponds to the jump process satisfying 
the non-resonance condition (see M. Arisawa and P.-L. Lions \cite{al}). 
At first sight, the formal asymptotic expansion (\ref{f1}) used for the $\a$-stable L{\'e}vy operator seems to be unapplicable for the measures in Examples 1-4. 
However, under the condition (A), by using the constant 
 $\a$ in it, we can still use the 
expansion (\ref{f1}) : 
$$
	u_{\e}(x)=\overline{u}(x)+  \e^{\a} v(\frac{x}{\e}) + o(\e^{\a})\qquad x\in {\bf R^N}. 
$$

We introduce the formal derivatives of $u_{\e}$ into (\ref{homo}) (resp. (\ref{homo2})). 
From  the condition (A)  (\ref{bound}),  (\ref{lim}),  we remark 
\begin{equation}\label{mc1}
	\lim_{\e\to 0} \int_{|z|\leq 1} |z|^{\gamma} \e^{M+\a} q(\e z) dz = \int_{|z|\leq 1} |z|^{\gamma} dq_0(z), 
\end{equation}
\begin{equation}\label{mc2}
	\lim_{\e\to 0} \int_{|z|> 1} |z|^{\gamma -1} \e^{M+\a} q(\e z) dz = \int_{|z| > 1} |z|^{\gamma -1} dq_0(z), 
\end{equation}
 for $\a\in (1,2)$ with $\gamma=2$ (resp. $\a\in (0,1)$ with $\gamma=1$). 
We get formally the following ergodic cell problem: for any fixed $x\in \Omega$ and for the given 
$$
	I_1= \int_{{\bf R^M}} 
	[\ou(x+\beta(z))-\ou(x)-   {\bf 1}_{|z|\leq 1} \la \n \ou(x),\beta(z) \ra] dq(z), 
$$ 
( resp. 
$$
	I_2= \int_{{\bf R^M}} 
	[\ou(x+\beta(z))-\ou(x)] dq(z), 
$$ 
) find a unique number 
$d_{I_1}$ (resp. $d_{I_2}$) such that the following problem has at least one periodic viscosity solution $v(y)$:
\begin{equation}\label{cell}
	d_{I_1}-a(y) \int_{{\bf R^M}} 
	[v(y+\beta(z))-v(y)- \la \n v(y),\beta(z) \ra] dq_0(z)  -a(y)I_1
\end{equation}
$$\qquad\qquad\qquad\qquad\qquad\qquad\qquad\qquad\qquad\qquad\qquad
	-f(y)=0\qquad \hbox{in}\quad {\bf T^N}, 
$$
(resp. 
\begin{equation}\label{cell2}
	d_{I_2}-a(y) \int_{{\bf R^M}} 
	[v(y+\beta(z))-v(y)] dq_0(z)  -a(y)I_2 -f(y)=0\qquad \hbox{in}\quad {\bf T^N}, 
\end{equation}
) 
provided that $dq_0(z)$ (the rescaled measure  defined in (\ref{q0})) satisfies  (\ref{radon}) with $\gamma=2$
 (resp. $\gamma=1$). 
 In some cases, we can only find the unique number $d_{I_1}$ (resp. $d_{I_2}$) which satisfies the following weaker property. For the case of (\ref{cell}), $d_{I_1}$ is the unique number such that for any $\delta>0$, there exist a subsolution $v_{\delta}$ and a supersolution $v^{\delta}$ of 
$$
	d_{I_1}-a(y) \int_{{\bf R^M}} 
	[v_{\delta}(y+\beta(z))-v_{\delta}(y)- \la \n v_{\delta}(y),\beta(z) \ra] dq_0(z)  -a(y)I_1
$$
$$\qquad\qquad\qquad\qquad\qquad\qquad\qquad\qquad\qquad\qquad\qquad
	-f(y)\leq \delta\qquad \hbox{in}\quad {\bf T^N},
$$
and
$$
	d_{I_1}-a(y) \int_{{\bf R^M}} 
	[v^{\delta}(y+\beta(z))-v^{\delta}(y)- \la \n v^{\delta}(y),\beta(z) \ra] dq_0(z)  -a(y)I_1
$$
$$\qquad\qquad\qquad\qquad\qquad\qquad\qquad\qquad\qquad\qquad\qquad
	-f(y)\geq -\delta\qquad \hbox{in}\quad {\bf T^N}. 
$$
The weaker version of (\ref{cell2}) will be stated later in \S 4. 
As remarked in \cite{al} for the case of partial differential equations, the existence of the unique number $d_{I_1}$ (resp. $d_{I_2}$) is shown by the strong maximum principle (SMP in short) for the L{\'e}vy operator. Since  the L{\'e}vy  density $dq_0(z)$ in (\ref{cell}) (resp.  (\ref{cell2})) is possibly degenerate,  we must establish a new SMP for our present purpose. We shall give a general sufficient condition for the SMP  in \S 2 (the condition (B)), in terms of 
 the controllability of the jump process: $x\to x+\beta(z)$ ($z\in S_0$). \\

$\quad$Although we have stated our problem in linear cases, for the reason of the simplicity, the present method is applicable to  nonlinear homogenization 
problems. \\

{\bf Example 5.} Let $\Omega\subset {\bf R^3}$, $\beta_1$$:{\bf R}\to {\bf R^3}$, $\beta_2$$:{\bf R^2}\to {\bf R^3}$ be such that 
$$
	\beta_1(z')=(0,0,z')\quad \forall z'\in {\bf R},\quad \beta_2(z'')=(z_1'',z_2'',0)\quad \forall z''=(z_1'',z_2'')\in {\bf R^2}. 
$$
Consider 
$$
	u_{\e}(x)+\max\{
	-a(\frac{x}{\e})\int_{{\bf R}}  [u_{\e}(x+\beta_1(z'))-u_{\e}(x)- {\bf 1}_{|z'|\leq 1}\la \n u_{\e}(x),\beta_1(z') \ra] dq_1(z'), 
$$
$$
	-a(\frac{x}{\e})\int_{{\bf R^2}}  [u_{\e}(x+\beta_2(z''))-u_{\e}(x)- {\bf 1}_{|z''|\leq 1}\la \n u_{\e}(x),\beta_2(z'') \ra] dq_2(z'') \}  
$$
\begin{equation}\label{ex4}
\qquad\qquad\qquad\qquad\qquad\qquad\qquad\qquad\qquad\qquad\quad
	-f(\frac{x}{\e})=0\quad \hbox{in}\quad \Omega, 
\end{equation} 
with the Dirichlet condition (\ref{bc}). 
Here, $dq_1(z')$, $dq_2(z'')$ are respectively a one-dimensional and a two-dimensional L{\'e}vy measures, and further detailed assumptions will be given later. We shall give the effective equation for this homogenization problem in \S 5. \\

The plan of this paper is the following. In \S 2, we state the SMP for  L{\'e}vy operators with degenerate densities satisfying a quite general 
condition (B) given in below. In \S 3, under the condition (B), we solve the ergodic cell problems (\ref{cell}) and (\ref{cell2}). In \S 4, the homogenization problem (\ref{homo}) and (\ref{homo2})  are solved rigorously. In \S 5, a generalization to nonlinear problems, such as Example 5, is indicated. In \S 6, the definitions of viscosity solutions for the integro-differential equations  with L{\'e}vy operators  are reviewed for the purpose of the readers. 
Throughout this paper, the notions of the subsolution and the supersolution mean  the viscosity subsolution and the viscosity supersolution, respectively. We denote by $USC({\bf R^N})$ and by $LSC({\bf R^N})$  the set of all upper semicontinuous functions on ${\bf R^N}$, and the set of all lower semicontinuous functions on ${\bf R^N}$, respectively. For $x\in {\bf R^N}$ we denote by $B_{r}(x)$ a ball centered at $x$ with radius $r>0$. \\

\section{Strong maximum principle in ${\bf T^N}$}.\\
 
$\quad$In this section, we establish the SMP for  L{\'e}vy operators with asymmetric, degenerate densities. We  use this result to solve the ergodic cell problem in \S 3. Our presentation is slightly more general than necessary. Let $H(y,p)$ be a continuous real valued function defined in ${\bf R^N}\times{\bf R^N}$, periodic in $y$ with the period ${\bf T^N}$, satisfying
\begin{equation}\label{H}
	H(y,0)\geq 0 \qquad \forall y\in {\bf T^N}. 
\end{equation}
We consider
$$
	H(y,\n u)-a(y)\int_{{\bf R^M}} 
	[u(y+\beta(z))-u(y)- \la \n u(y),\beta(z) \ra] dq_0(z)  =0
$$
\begin{equation}\label{smp}
\qquad\qquad\qquad\qquad\qquad\qquad\qquad\qquad\qquad\qquad\qquad\qquad
	 \hbox{in}\quad {\bf T^N}, 
\end{equation}
and 
\begin{equation}\label{smp2}
	H(y,\n u)-a(y)\int_{{\bf R^M}} 
	[u(y+\beta(z))-u(y)] dq_0(z)  =0\quad
	 \hbox{in}\quad {\bf T^N}, 
\end{equation}
where $\beta(z)$ satisfies (\ref{beta}),  $a(y)$ satisfies (\ref{holder1}), and $dq_0(z)$ satisfies (\ref{radon}) with $\gamma=2$ in the case of (\ref{smp}), with $\gamma=1$ in the case of (\ref{smp2}) respectively. 
We assume the following  condition. \\

(B) For any two points $y$, $y' \in {\bf T^N}$, there exist a finite number of points $y_1$, ..., $y_m\in {\bf T^N}$ such that $y_1=y$, $y_m=y'$, and for any 
$m$ positive numbers $\e_i>0$ $(1\leq i\leq m)$, we can take 
subsets $J_i\subset S_0$$=\hbox{supp}(dq_0(z))$ ($1\leq \forall i\leq m-1$) satisfying 
\begin{equation}\label{C}
	y_i+\beta(z)\in B_{\e_i}(y_{i+1}) \quad \forall z\in J_i;\quad \int_{J_i} 1 dq_0(z)> 0 \quad 1\leq \forall i\leq m. 
\end{equation}

The condition (B) describes the controllability of the jump process $y$$\to$$y+\beta(z)$ ($z\in S_0$). \\

{\bf Theorem 2.1.$\quad$}
\begin{theorem}  Let $u\in USC({\bf R^N})$  be a viscosity subsolution of (\ref{smp}) (resp. (\ref{smp2})). Assume that  (\ref{beta}), (\ref{holder1}), (\ref{H}) hold, and that $dq_0(z)$ satisfies the condition (B) and  (\ref{radon}) with $\gamma=2$ (resp. $\gamma=1$). 
  If $u$ attains a maximum at $\oy$ in ${\bf T^N}$, then $u$ is constant in ${\bf T^N}$. 
\end{theorem} 

$Proof.$ Let $u(\oy)=M$, and put $\Omega_0=$$\{y\in {\bf T^N}|\quad u\equiv M\}$. Assume that $\Omega_0^c\neq \emptyset$, and we shall lead a contradiction. 
Take a point $y'\in \Omega_0^c$, and remark that $u(y')<M$. From the condition (B), we can take a finite number of points, $y_1$,..., $y_m$$\in {\bf T^N}$  such that $y_1=\oy$, $y_m=y'$, $m$ positive numbers $\e_i$ $(1\leq i\leq m)$, and $m-1$ subsets $J_i\subset S_0$ which satisfy (\ref{C}). There exists a number $k$  such that $1\leq k<m$, with which $y_k\in \Omega_0$ and $y_{k+1}\in \Omega_0^c$. 
Since $\Omega_0^c$ is open, we can take $\e_{k}>0$ small enough so that $B_{\e_k}(y_{k+1})\subset \Omega_0^c$. From the condition (B), there exists $J_k\subset S_0=\hbox{supp}(dq_0(z))$ such that $\int_{J_k} 1 dq_0(z) >0$, and 
$$
	y_k+\beta(z)\in U_{\e_k}(y_{k+1}) \quad \forall z\in J_k. 
$$ 
Thus, we can take $\delta_k>0$ such that 
\begin{equation}\label{dk}
	u(y_k+\beta(z))< M-\delta_k \quad \forall z\in J_k.
\end{equation}
For the constant function $\phi(y)\equiv M$ ($y\in {\bf T^N}$), since $u-\phi$ takes a maximum at $y_k$, from the definition of the viscosity subsolution (see Definition C in \S 6), by using $\n \phi(y_k)=0$, 
 we have
$$
	H(y_k,0)-a(y_k)\int_{{\bf R^M}} 
	[u(y_k+\beta(z))-u(y_k)- \la 0,\beta(z) \ra] dq_0(z)  \leq 0. 
$$
( resp. 
$$
	H(y_k,0)-a(y_k)\int_{{\bf R^M}} 
	[u(y_k+\beta(z))-u(y_k)] dq_0(z)  \leq 0.)
$$
 From (\ref{holder1}), (\ref{H}), and from the fact that $u(y_k)=M > u(y_k+\beta(z))$ for any $z\in \hbox{supp}(dq_0(z))$, the above leads to 
$$
	-\int_{J_k} 
	[u(y_k+\beta(z))-u(y_k)] dq_0(z)  \leq 0. 
$$
However, from the condition (B), this contradicts to (\ref{dk}), since $-\int_{J_k} [u(y_k+\beta(z))-u(y_k)] dq_0(z)$$\geq \delta_k\int_{J_k}1dq_0(z)>0$. 
Therefore, $\Omega_0^c=\emptyset$ must hold. \\

{\bf Remarks 2.1.} 1. Consider the jump process: $y\to y+\beta(z)$ ($z\in S_0=\hbox{supp}(dq_0(z))$) in ${\bf T^N}$, where $dq_0(z)$ is either one of the measures 
defined in Examples 1-4. Then, it is easy to see that the condition (B) is satisfied by each of the measures $dq_0(z)$. (Remark that in Example 3, for $y\in {\bf T^2}$ fixed, the set $\{y+(z,\gamma z)| z\in {\bf R=S_0}\}$ is dense in ${\bf T^2}$ for $\gamma>0$ is irrational.)\\

2. Let $M=N$, and $\beta(z)=z$. If for some $r>0$, $B_r(0)\subset dq_0(z)$, then the condition (B) is satisfied. \\

3. The SMP in Theorem 2.1 can be stated in parallel  for a supersolution $u\in LSC({\bf R^N})$ of  (\ref{smp}) (resp. (\ref{smp2})), i.e. if $u$ attains a minimum 
 at $\oy\in {\bf T^N}$, then $u$ is a constant function.\\

4. Let us replace the L{\'e}vy operator in (\ref{smp}) to the following : 
$$
	\int_{{\bf R^M}} [u(y+\beta(z))-u(y)- {\bf 1}_{|z|\leq 1} \la \n u(y),\beta(z) \ra] dq_0(z),
$$
where $dq_0(z)$ satisfies (\ref{radon}) with $\gamma=2$. Then the SMP holds for the above operator, under the condition (B), too. \\

\section{Ergodic problem.}
In this section, we study the ergodic problem of the jump process: $x \to x+\beta(z)$ 
$(z\in {supp}(dq_0(z)))$.   
For $\l>0$, we consider 
\begin{equation}\label{ergodic}
	\l v_{\l}(y)-a(y)\int_{{\bf R^M}}  [v_{\l}(y+\beta(z))-v_{\l}(y)- \la \n v_{\l}(y),\beta(z) \ra] dq_0(z) 
\end{equation}
$$
	\qquad\quad\qquad\qquad\quad\qquad\qquad\quad\qquad\qquad\quad\qquad
	-f_0(y)=0\quad \hbox{in}\quad {\bf T^N}. 
$$
( resp. 
\begin{equation}\label{ergodic2}
	\l v_{\l}(y)-a(y)\int_{{\bf R^M}}  [v_{\l}(y+\beta(z))-v_{\l}(y)] dq_0(z) -f_0(y)=0\quad \hbox{in}\quad {\bf T^N}. 
\end{equation}
) It is known that there exists a unique periodic viscosity solution $v_{\l}$ of (\ref{ergodic}) (resp. (\ref{ergodic2}))
(see \cite{arnewdef}, \cite{arcorrig}, and \cite{bi}). \\

{\bf Theorem 3.1.$\quad$} 
\begin{theorem}  Let $v_{\l}$  be a viscosity solution of (\ref{ergodic}) (resp. (\ref{ergodic2})). Assume that (\ref{beta}), (\ref{holder1}) hold, that $f_0$ satisfies (\ref{holder2}), that $dq_0(z)$ satisfies the condition (B) and (\ref{radon}) with $\gamma=2$ (resp. $\gamma=1$).  
  Then, there exists a unique real number $d$ such that 
\begin{equation}\label{limit}
	\lim_{\l\to 0}\l v_{\l}(y)=d\qquad \hbox{uniformly in}\quad  {\bf T^N}. 
\end{equation}
The number $d$ is characterized by the following property: for any $\delta>0$ there exists a subsolution $v_{\delta}$ and a supersolution $v^{\delta}$ of 
\begin{equation}\label{aergodic1}
	d-a(y)\int_{{\bf R^M}}  [v_{\delta}(y+\beta(z))-v_{\delta}(y)- \la \n v_{\delta}(y),\beta(z) \ra] dq_0(z) -f_0(y)\leq \delta, 
\end{equation}
\begin{equation}\label{aergodic2}
	d-a(y)\int_{{\bf R^M}}  [v^{\delta}(y+\beta(z))-v^{\delta}(y)- \la \n v^{\delta}(y),\beta(z) \ra] dq_0(z) -f_0(y)\geq -\delta, 
\end{equation}
( resp. 
\begin{equation}\label{aergodic12}
	d-a(y)\int_{{\bf R^M}}  [v_{\delta}(y+\beta(z))-v_{\delta}(y)] dq_0(z) -f_0(y)\leq \delta, 
\end{equation}
\begin{equation}\label{aergodic22}
	d-a(y)\int_{{\bf R^M}}  [v^{\delta}(y+\beta(z))-v^{\delta}(y)] dq_0(z) -f_0(y)\geq -\delta, 
\end{equation}
) in ${\bf T^N}$ respectively. 
\end{theorem} 

$Proof.$ We prove (\ref{limit}) for the problem (\ref{ergodic}). The proof for (\ref{ergodic2}) is similar and we do not write it here. 
 We multiply (\ref{ergodic}) by $\l>0$, and put $m_{\l}=\l v_{\l}$. We have
\begin{equation}\label{ml}
	\l m_{\l}(y)-a(y)\int_{{\bf R^M}}  [m_{\l}(y+\beta(z))-m_{\l}(y)- \la \n m_{\l}(y),\beta(z) \ra] dq_0(z)
\end{equation}
$$
	\qquad\quad\qquad\qquad\quad\qquad\qquad\quad\qquad\qquad\quad\qquad
	 -\l f_0(y)=0\quad \hbox{in}\quad {\bf T^N}. 
$$
We claim that the following holds. \\

{\bf Lemma 3.2.$\quad$}
\begin{theorem}  Let the assumptions in Theorem 3.1 hold. \\
(i)  There exists a constant $M>0$ such that the following hold: 
\begin{equation}\label{est11}
	 |m_{\l}|_{L^{\infty}}\leq M \quad \forall \l\in (0,1). 
\end{equation}
(ii) For any $\theta\in (0,\min\{\theta_1,\theta_2\})$, there exists a constant $C_{\theta}>0$ such that 
\begin{equation}\label{est22}
	|m_{\l}(y)- m_{\l}(y')| \leq C_{\theta}|y-y'|^{\theta} \quad \forall y,y'\in {\bf T^N}, \quad  \forall \l\in (0,1). 
\end{equation}
The constants $M$, $C_{\theta}>0$ are independent on $\l\in (0,1)$. 
\end{theorem} 

We admit the above estimates for a while, which we shall prove later. By Lemma 3.2 ($m_{\l}=\l v_{\l}$), from the Ascoli-Arzera lemma we can take a sequence $\l'\to 0$ such that 
$$
	\l' v_{\l'}(y) \to \exists d(y) \quad \hbox{as}\quad \l' \to 0,\quad \hbox{uniformly in}\quad {\bf T^N},
$$
where $d(y)$ is a H{\"o}lder continuous, periodic function satisfying (\ref{est22}). To see that $d(y)$ is constant, we multiply (\ref{ergodic}) by $\l'>0$, and  tend $\l'$ to zero. By using (\ref{est11}), from the stability of viscosity solutions, we get 
$$
	-\int_{{\bf R^M}}  [d(y+\beta(z))-d(y)- \la \n d(y),\beta(z) \ra] dq_0(z) \leq 0 \quad \hbox{in}\quad {\bf T^N}. 
$$
Hence, from the SMP in Theorem 2.1, $d(y)$ is constant, i.e. $d(y)\equiv d$ 
for some real number $d$. Next, assume that there exists another sequence $\l''\to 0$ and another number $d'$ such that 
$$
	\l'' v_{\l''}(y) \to d' \quad \hbox{as}\quad \l'' \to 0,\quad \hbox{uniformly in}\quad {\bf T^N}. 
$$
Without loss of generality, we may assume that $d'< d$. 
For arbitrary small $\mu>0$, by taking $\l'>0$ and $\l''>0$ small enough, we have the following two inequalities. 
$$
	d-a(y)\int_{{\bf R^M}}  [v_{\l'}(y+\beta(z))-v_{\l'}(y)- \la \n v_{\l'}(y),\beta(z) \ra] dq_0(z)-f_0(y)\leq \frac{\mu}{2},
$$
$$
	d'-a(y)\int_{{\bf R^M}}  [v_{\l''}(y+\beta(z))-v_{\l''}(y)- \la \n v_{\l''}(y),\beta(z) \ra] dq_0(z) -f_0(y)\geq -\frac{\mu}{2}.
$$
We shall write $\underline{w}=v_{\l'}$, $\overline{w}=v_{\l''}$. By adding a constant if necessary,  
 we may assume that 
\begin{equation}\label{bigger}
	\underline{w}(y) >  \overline{w}(y) \quad  \forall y\in {\bf T^N}. 
\end{equation}
We take $\l>0$ small enough  so that $|\l \underline{w}|_{L^{\infty}}$, $|\l \overline{w}|_{L^{\infty}}$$<\frac{\mu}{2}$.  Then, 
 $\underline{w}$ and $\overline{w}$ satisfy respectively 
$$
	\l \underline{w}(y)-a(y)\int_{{\bf R^M}}  [\underline{w}(y+\beta(z))-\underline{w}(y)- \la \n \underline{w}(y),\beta(z) \ra] dq_0(z)
$$
$$\qquad\quad\qquad\qquad\quad\qquad\qquad\quad\qquad\qquad\quad\qquad\quad\qquad
	 +d -f_0(y)\leq {\mu},
$$
$$
	\l \overline{w}(y)-a(y)\int_{{\bf R^M}}  [\overline{w}(y+\beta(z))-\overline{w}(y)- \la \n \overline{w}(y),\beta(z) \ra] dq_0(z) 
$$
$$\qquad\quad\qquad\qquad\quad\qquad\qquad\quad\qquad\qquad\quad\qquad\quad\qquad
	 + d' -f_0(y)\geq -{\mu}.
$$
From the comparison principle (see \cite{arnewdef}, \cite{arcorrig}, \cite{bi}), we get 
$$
	\l (\underline{w}(y)-   \overline{w}(y)) \leq d'-d+ 2\mu\quad \forall y\in {\bf T^N}, 
$$
which contradicts to (\ref{bigger}), for $\mu>0$ small enough. Therefore, $d=d'$ should hold, and we proved the claim. \\

$Proof$ $of$ $Lemma $ $3.2.$
(i) The uniform bound for $|m_{\l}|_{L^{\infty}}$ ($\forall \l\in (0,1)$) is clear from the comparison principle for (\ref{ml}),
 i.e. $|\l m_{\l}|_{L^{\infty}} \leq |\l f_0 |_{L^{\infty}}$. \\

(ii) We show the inequality by the contradiction argument. Let $r>0$ be a fixed number to be determined later. Put 
\begin{equation}\label{ctheta}
	C_{\theta}=\frac{2M}{r^{\theta}}. 
\end{equation}
Assume that there exist $\oy$, $\oy'\in {\bf T^N}$ such that 
\begin{equation}\label{assump}
	|m_{\l}(\oy)-m_{\l}(\oy')|>C_{\theta}|\oy-\oy'|^{\theta}, 
\end{equation}
and we shall lead a contradiction. Remark that $|\oy-\oy'|<r$ must hold. Put 
$$
	\Phi(y,y')=m_{\l}(y)-m_{\l}(y')-C_{\theta}|y-y'|^{\theta} \quad y,y'\in {\bf T^N}. 
$$
Let $(\hy,\hy')$ be a maximum point of $\Phi$ in ${\bf T^N}$. We may assume that $\Phi(\hy,\hy')$ is the strict maximum.
 Put $\phi(y,y')=C_{\theta}|y-y'|^{\theta}$, $p=\n_y\phi(\hy,\hy')$, $Q=\n^2_y \phi(\hy,\hy')$. From the definition of the viscosity solution, we get 
$$
	\l m_{\l}(\hy)-a(\hy) \int_{{\bf R^M}}  [m_{\l}(\hy+\beta(z))-m_{\l}(\hy)- \la p,\beta(z) \ra] dq_0(z) \leq \l f_0(\hy), 
$$
$$
	\l m_{\l}(\hy')-a(\hy')\int_{{\bf R^M}}  [m_{\l}(\hy'+\beta(z))-m_{\l}(\hy')- \la p,\beta(z) \ra] dq_0(z) \geq \l f_0(\hy'). 
$$

By deviding the above two inequalities  by $a(\hy)$ and $a(\hy')$ respectively, and then by taking the difference of them, we have 
 $$
	\frac{\l m_{\l}(\hy)}{a(\hy)}-\frac{\l m_{\l}(\hy')}{a(\hy')}- 
	\int_{{\bf R^M}}  [m_{\l}(\hy+\beta(z))-m_{\l}(\hy) \quad\qquad\qquad\qquad
$$
$$\qquad\quad\qquad\qquad\quad
	- m_{\l}(\hy'+\beta(z))+ m_{\l}(\hy')] dq_0(z) 
	\leq \frac{\l f_0(\hy)}{a(\hy)} - \frac{\l f_0(\hy')}{ a(\hy')}.  
$$
Since for any $z\in {\bf R^M}$, 
$$
	m_{\l}(\hy)-m_{\l}(\hy')-C_{\theta}|\hy-\hy'|^{\theta} \geq m_{\l}(\hy+\beta(z))-m_{\l}(\hy'+\beta(z))-C_{\theta}|\hy-\hy'|^{\theta}, 
$$
the preceding inequality leads to

$$
	\l a(\hy')m_{\l}(\hy)-\l a(\hy)m_{\l}(\hy')\leq \l a(\hy')f_0(\hy)- \l a(\hy)f_0(\hy'), 
$$
which leads to
$$
	a(\hy')(m_{\l}(\hy)-m_{\l}(\hy')) \qquad\quad\qquad\qquad\quad\qquad\quad\qquad\qquad\quad\qquad\quad\
$$
$$\qquad\quad\qquad\qquad
	\leq (a(\hy)-a(\hy')) m_{\l}(\hy') + a(\hy')(   f_0(\hy)- f_0(\hy'))+(a(\hy')-a(\hy)) f_0(\hy'). 
$$
Thus, from (\ref{holder1}), (\ref{holder2}),  (\ref{assump}), since $(\hy,\hy')$ is the maximum point of $\Phi$, the above leads to 
$$
	C_{\theta}|\hy-\hy'|^{\theta} \leq L'(|\hy-\hy'|^{\theta_1}+ |\hy-\hy'|^{\theta_2}), 
$$
where $L'=a_0^{-1}L(M+||a||_{L^{\infty}({\bf T^N})}$$+ ||f_0||_{L^{\infty}({\bf T^N})})$. 
Therefore,  from (\ref{ctheta}), since $\theta\in$$(0,\min\{\theta_1,\theta_2\})$ and $|\hx-\hx'|<r$, 
$$
	2M\leq L'(|\hy-\hy'|^{\theta_1-\theta}r^{\theta}+|\hy-\hy'|^{\theta_2-\theta}r^{\theta})\leq L'(r^{\theta_1}+r^{\theta_2}). 
$$ By taking $r>0$ small enough so that $r^{\theta_1}+ r^{\theta_2}<2ML'^{-1}$, we get a desired 
contradiction. This shows the existence of $C_{\theta}>0$ such that (ii) holds. Moreover, the constant $C_{\theta}$ does not depend on $\l\in (0,1)$.\\

{\bf Corollary 3.3.$\quad$}
\begin{theorem}  (i) Let $v_{\l}$ be the solution of (\ref{ergodic}) with $dq_0(z)$ and $\beta(z)$ given either one of the following : Example 1 with $\a\in (1,2)$, 
 Examples 2 and 3 with $\a\in (1,2)$, and Exmple 4 with $\a\in (1,2)$. Then, there exists a unique constant $d$ 
such that (\ref{limit}) holds. \\
(ii) Let $v_{\l}$ be the solution of (\ref{ergodic2}) with $dq_0(z)$ and $\beta(z)$ given either one of the following : Example 1 with $\a\in (0,1)$, 
 Example 3 with $\a\in (0,1)$, and Exmple 4 with $\a\in (0,1)$. Then, there exists a unique constant $d$ 
such that (\ref{limit}) holds. \\
\end{theorem} 
$Proof.$ As we have seen in Remarks 2.1, each of the measures $dq_0(z)$ in Examples 1-4 satisfies the condition (B). Hence, the claim follows from Theorem 3.1.\\

{\bf Remarks 3.1.} 1. The SMP (Theorem 2.1)  is  essential to prove the existence of the ergodic number $d$ in Theorem 3.1. 
 \\

2. We can generalize Theorem 3.1, by adding a fully nonlinear degenerate elliptic second-order operator $F(x,\n u,\n^2 u)$ to (\ref{ergodic}) (resp. (\ref{ergodic2})). \\

\section{Homogenizations.}

$\quad$ In this section, we give our main results of the homogenization problems (\ref{homo})-(\ref{bc}) and (\ref{homo2})-(\ref{bc}) in Theorems 4.6 and 4.8 respectively.   
Throughout this section, we assume that the condition (A) holds.  
Let $u_{\e}$ be the solution of (\ref{homo})-(\ref{bc})  (resp. (\ref{homo2})-(\ref{bc}) ). By introducing the formal asymptotic expansion (\ref{f1}): 
$$
	u_{\e}(x)=\overline{u}(x)+  \e^{\a} v(\frac{x}{\e}) + o(\e^{\a})\qquad x\in {\bf R^N},
$$
into (\ref{homo}) (resp. (\ref{homo2})), by using the homogeneity of $\beta$ in (\ref{beta}), by remarking that (\ref{mc1}) and (\ref{mc2}) hold, 
we get the following cell problem (\ref{cell}):  
$$
	d_{I_1}-a(y)\int_{{\bf R^M}} 
	[v(y+\beta(z))-v(y)- \la \n v(y),\beta(z) \ra] dq_0(z)\qquad\qquad\qquad
$$
$$
	\qquad\qquad\qquad\qquad\qquad\qquad
	 -a(y)I_1-f(y)=0\qquad \hbox{in}\quad {\bf T^N}, 
$$
where 
$$
	I_1=I_1[\ou](x)= \int_{{\bf R^M}} 
	[\ou(x+\beta(z))-\ou(x)- {\bf 1}_{|z|\leq 1}\la \n \ou(x),\beta(z) \ra] dq(z), 
$$
( resp. (\ref{cell2}):  
$$
	d_{I_2}-a(y)\int_{{\bf R^M}} 
	[v(y+\beta(z))-v(y)] dq_0(z)-a(y)I_2-f(y)=0\qquad \hbox{in}\quad {\bf T^N}, 
$$
where 
$$
	I_2=I_2[\ou](x)= \int_{{\bf R^M}} [\ou(x+\beta(z))-\ou(x)] dq(z), 
$$
) provided that $dq_0(z)$ satisfies (\ref{radon}) with $\gamma=2$ (resp. $\gamma=1$). 
 Remark that according to the condition (A), the L{\'e}vy measure $dq(z)$ in (\ref{homo}) (resp. (\ref{homo2})) is transformed to $dq_0(z)$ in the cell problem (\ref{cell}) (resp. (\ref{cell2})).  
For any $I_1\in {\bf R}$ (resp. $I_2\in {\bf R}$), from Theorem 3.1 (with $f_0(y)=$$a(y)I_i +f(y)$ (i=1,2)), there exists a unique number $d_{I_1}$ (resp. $d_{I_2}$) such that for any $\delta>0$ there exist $v_{\delta}$ a  periodic subsolution and $v^{\delta}$ a periodic supersolution of 
$$
	d_{I_1}-a(y)\int_{{\bf R^M}} 
	[v_{\delta}(y+\beta(z))-v_{\delta}(y)- \la \n v_{\delta}(y),\beta(z) \ra] dq_0(z) 
$$
$$\qquad\qquad\qquad\qquad\qquad\qquad\qquad
	 -a(y)I_1-f(y)\leq \frac{\delta}{2} \qquad \hbox{in}\quad {\bf T^N}, 
$$
$$
	d_{I_1}-a(y)\int_{{\bf R^M}} 
	[v^{\delta}(y+\beta(z))-v^{\delta}(y)  - \la \n v_{\delta}(y),\beta(z) \ra  ] dq_0(z) 
$$
$$\qquad\qquad\qquad\qquad\qquad\qquad\qquad
	 -a(y)I_1 -f(y)\geq - \frac{\delta}{2} \qquad \hbox{in}\quad {\bf T^N}. 
$$
( resp. 
$$
	d_{I_2}-a(y)\int_{{\bf R^M}} 
	[v_{\delta}(y+\beta(z))-v_{\delta}(y)] dq_0(z) -a(y)I_2-f(y)\leq \frac{\delta}{2} \qquad \hbox{in}\quad {\bf T^N}, 
$$
$$
	d_{I_2}-a(y)\int_{{\bf R^M}} 
	[v^{\delta}(y+\beta(z))-v^{\delta}(y)] dq_0(z) -a(y)I_2 -f(y)\geq - \frac{\delta}{2} \qquad \hbox{in}\quad {\bf T^N}. 
$$
) For the later purpose, let us regularize $v_{\delta}$ and $v^{\delta}$: for $\nu>0$, define 
$$
	v_{\delta}^{\nu}(x)=\sup_{|y-x|\leq \nu}\{v_{\delta}(y)-\frac{1}{\nu^2}|y-x|^2\} \quad (\hbox{sup convolution}),
$$
$$
	v_{\nu}^{\delta}(x)=\inf_{|y-x|\leq \nu}\{v^{\delta}(y)+\frac{1}{\nu^2}|y-x|^2\}\quad (\hbox{inf convolution}). 
$$ 
Put 
$\uv=v_{\delta}^{\nu}$, $\ov=v_{\nu}^{\delta}$.  
It is known that $\uv$ is semiconvex, $\ov$ is semiconcave, and both are Lipschitz continuous (see \cite{users}, W.H. Fleming and H.M. Soner \cite{fs}). Moreover, since $\lim_{\nu\to 0}v_{\delta}^{\nu}=v_{\delta}$, $\lim_{\nu\to 0}v^{\delta}_{\nu}=v^{\delta}$ uniformly in ${\bf T^N}$, 
for any $\delta>0$, we can take 
$\nu>0$ such that  $\uv$  and $\ov$ are respectively  a subsolution and a supersolution of the following : 
\begin{equation}\label{approxcell1}
	d_{I_1}-a(y)\int_{{\bf R^M}} 
	[\uv(y+\beta(z))-\uv(y)- \la \n \uv(y),\beta(z) \ra] dq_0(z) 
\end{equation}
$$\qquad\qquad\qquad\qquad\qquad\qquad\qquad
	 -a(y)I_1-f(y)\leq \delta \qquad \hbox{in}\quad {\bf T^N}, 
$$
\begin{equation}\label{approxcell2}
	d_{I_1}-a(y)\int_{{\bf R^M}} 
	[\ov(y+\beta(z))-\ov(y)- \la \n \ov(y),\beta(z) \ra] dq_0(z) 
\end{equation}
$$\qquad\qquad\qquad\qquad\qquad\qquad\qquad
	-a(y)I_1-f(y)\geq -\delta \qquad \hbox{in}\quad {\bf T^N}, 
$$
( resp. 
\begin{equation}\label{approxcell12}
	d_{I_2}-a(y)\int_{{\bf R^M}} 
	[\uv(y+\beta(z))-\uv(y)] dq_0(z) -a(y)I_2-f(y)\leq \delta \quad \hbox{in}\quad {\bf T^N}, 
\end{equation}
\begin{equation}\label{approxcell22}
	d_{I_2}-a(y)\int_{{\bf R^M}} 
	[\ov(y+\beta(z))-\ov(y)] dq_0(z) -a(y)I_2-f(y)\geq -\delta \quad \hbox{in}\quad {\bf T^N}, 
\end{equation}
) (see for instance \cite{arcorrig}, \cite{fs}). We use the above approximated cell problem in place of (\ref{cell}) in the following argument. Define 
\begin{equation}\label{oI}
	\oI_1(I_1)= -d_{I_1} \quad \forall I_1\in {\bf R}\quad (resp. \quad \oI_2(I_2)= -d_{I_2} \quad \forall I_2\in {\bf R}), 
\end{equation}
where the right hand side is the unique number such that for any $\delta>0$, (\ref{approxcell1}) and (\ref{approxcell2}) (resp. (\ref{approxcell12}) and (\ref{approxcell22})) have a subsolution and a supersolution respectively. We prepare some lemmas which we use later. \\

{\bf Lemma 4.1.$\quad$ (\cite{arperiodic})}
\begin{theorem}  Assume that (\ref{beta}), (\ref{holder1}), (\ref{holder2}) and that $dq_0(z)$ satisfies the condition  (B) and  (\ref{radon}) with $\gamma=2$  (resp.    $\gamma=1$). Then, the function $\oI_1$ (resp. $\oI_2$) defined in (\ref{oI}) is continuous and satisfies  the following property: there exists $\Theta>0$ such that 
\begin{equation}\label{subelliptic}
	\oI_1(I+I') - \oI_1(I) \leq  - \Theta I' 
	\quad (resp. \oI_2(I+I') - \oI_2(I) \leq  - \Theta I' )
	\quad\forall I\in {\bf R},\forall I'\geq 0. 
\end{equation}
\end{theorem}

The above result was presented in \cite{arperiodic}, which was originally given in \cite{ev1} for the PDE case. 
The proof does not differ so much from \cite{arperiodic}, \cite{ev1}, and we omit it here. \\

{\bf Remarks 4.1.} Let $u\in C^2(\bf R^N)$. Then, by putting 
$$
	I_1=I_1[u](x)= \int_{{\bf R^M}} 
	[u(x+\beta(z))-u(x)- {\bf 1}_{|z|\leq 1} \la \n u(x),\beta(z) \ra] dq(z)
$$
( resp. 
$$
	I_2=I_2[u](x)= \int_{{\bf R^M}} [u(x+\beta(z))-u(x)] dq(z)
$$
) into $\oI_1$ (resp. $\oI_2$), the map:
$u\to \oI_1(I_1[u](x))$ (resp. $u\to \oI_2(I_2[u](x))$) can be regarded as an integro-differential operator. The property (\ref{subelliptic}) implies that 
 $\oI_1(I_1[u](x))$  (resp.  $\oI_2(I_2[u](x))$) is subelliptic (see \cite{arperiodic}).\\

{\bf Lemma 4.2.$\quad$}
\begin{theorem}  Let $\oI_1$ (resp. $\oI_2$) be the functions defined in (\ref{oI}). Consider 
\begin{equation}\label{unique}
	u+\oI_1(I_1[u](x))=0 \quad(resp. u+\oI_2(I_2[u](x))=0) \quad \hbox{in}\quad \Omega,
\end{equation}
with (\ref{bc}). Let $u$, $v$ be respectively a subsolution and a supersolution of (\ref{unique})-(\ref{bc}). Then, $u\leq v$ in $\Omega$. 
Moreover, there exists a unique viscosity solution $u$  of (\ref{unique})-(\ref{bc}). 
\end{theorem}
$Proof.$ The comparison principle can be shown by the usual contradintion argument, from the subellipticity  (\ref{subelliptic}). 
The existence of the solution can be obtained by the Perron's method. This argument was done in 
  \cite{arnewdef}, \cite{arcorrig}, \cite{arperiodic}, \cite{bi} and we do not repeat it here.\\ 

We remind the following result in the convex analysis, which we cite without proof (see \cite{users}, \cite{fs} for details). For an upper or a lower semicontinuous function
 $\Phi$ defined in an open subset $\mathcal{O}$ in ${\bf R^n}$,  for $\rho>0$, put
$$
	M_{\rho}=\{\ox\in\mathcal{O} | \quad \exists p\in{\bf R^n}\quad \hbox{s.t.}\quad |p|\leq \rho, \quad \Phi(x)\leq \Phi(\ox)+\la p,x-\ox\ra\quad \forall x\in  \mathcal{O}\}. 
$$

{\bf Lemma 4.3.$\quad$ (\cite{users},\quad \cite{fs})}
\begin{theorem}  Let $\Phi$ be a semiconvex function in an open domain $\mathcal{O}$, and let $x'$ be a maximizer of $\Phi$ in $\mathcal{O}$ such that 
$$
	\mu=\sup_{\mathcal{O}} \Phi(x) -\sup_{\p \mathcal{O}} \Phi(x)= \Phi(x')-\sup_{\p \mathcal{O}} \Phi(x) >0.
$$
Then, the following hold.\\
(i) $\Phi$ is differentiable at $x'$ and $\n \Phi(x')$=0.\\
(ii) For any $m\in {\bf N}$, there exists $x_m\in M_{\frac{1}{m}}$ such that $\Phi$ is twice differentiable at $x_m$, $\lim_{m\to \infty} x_m=x'$, $\n^2 \Phi(x_m)\leq O$, $|\n \Phi(x_m)|\leq \frac{1}{m}$. For $p_m$$=\n \Phi(x_m)$, the function
$$
	\Phi_m(x)=\Phi(x)-\la p_m,x \ra
$$
takes a maximum at $x=x_m$. 
\end{theorem} 

{\bf Lemma 4.4.$\quad$} \begin{theorem}  Let $\uv(y)$ be a periodic semiconcave function defined in ${\bf T^N}$. Assume that for a function $\Psi(x)$$\in C^2({\bf R^N})$,  $\Psi(x)+\e^{\a}\uv(\frac{x}{\e})$ takes a global minimum at $\ox$.  
 Then, the following hold for any $z\in {\bf R^M}$, 
with a constant $C>0$ independent on $\e>0$ and $\ox$. \\
(i)
\begin{equation}\label{estv}
	-C\e^{2-\a} |z|^2 \leq 
	\ov(\frac{\ox}{\e}+\beta(z))-\ov(\frac{\ox}{\e})-\la \n_y \ov(\frac{\ox}{\e}), \beta(z) \ra
	\leq C|z|^2. 
\end{equation}
(ii)
 \begin{equation}\label{estv2}
	-C\e^{1-\a}|z| \leq 
	\ov(\frac{\ox}{\e}+\beta(z))-\ov(\frac{\ox}{\e})
	\leq C|z|. 
\end{equation}
\end{theorem} 

 $Proof$ $of$ $Lemma$ $4.4.$ (i) The second inequality comes from the semiconcavity of $\ov$ and (\ref{beta}). The first inequality is derived from the fact that 
$\Psi(x)+\e^{\a}\ov(\frac{x}{\e})$ takes a  global minimum at $\ox$. In fact, since $\Psi(x)+\e^{\a}\ov (\frac{x}{\e})$ is semiconcave, it is differentiable at $\ox$ and 
$\n \Psi(\ox)+$$\e^{-1+\a} \n_y \ov(\frac{\ox}{\e})$$=0$, 
$$
	\Psi(\ox)+\e^{\a}   \ov(\frac{\ox}{\e})
	\leq \Psi(\ox+\e\beta(z))+\e^{\a} \ov(\frac{\ox}{\e}+\beta(z))\quad \forall z\in {\bf R^M},
$$
 for any $\e>0$. Thus, we get 
$$
	\e^{\a}( \ov(\frac{\ox}{\e}+\beta(z)) - \ov(\frac{\ox}{\e}) -\la \e^{-1}\n_y \ov(\frac{\ox}{\e}),\e\beta(z)\ra)
	\qquad\qquad\qquad\qquad\qquad\qquad\qquad\qquad\qquad
$$
$$
	\geq -(\Psi(\ox+ \e\beta(z))-\Psi(\ox)-\la \n\Psi(\ox),\e\beta(z) \ra) \geq -\e^2|\beta(z)|^2 |\n^2\Psi(\ox+\mu\e\beta(z))|,
$$
where $\mu\in(0,1)$. 
From (\ref{beta}), the first inequality holds with a constant $C>0$ independent on $\e>0$ and $\ox$. \\
(ii) The second inequality comes from the Lipschitz continuity of $\ov$ and (\ref{beta}).  The first inequality is proved similarly to (i).\\

{\bf Lemma 4.5.$\quad$} \begin{theorem}  Let $\ov(y)$ be a periodic semiconcave function defined in ${\bf T^N}$. 
  Assume that for a function $\Psi(x)$$\in C^2({\bf R^N})$,  $\Psi(x)+g(\e)\ov(\frac{x}{\e})$ takes a  minimum at $\ox$.  
 Then, the following hold. \\
(i)  If $dq_0(z)$ satisfies (\ref{radon}) with $\gamma=2$,  
$$
	 \e^{\a}\int_{{\bf R^M}} [\ov(\frac{\ox+\beta(z)}{\e})-\ov(\frac{\ox}{\e})-{\bf 1}_{|z|\leq 1}\la \e^{-1}\n_y \ov(\frac{\ox}{\e}), \beta(z) \ra ]dq(z)
$$
\begin{equation}\label{est1}
	=  \int_{{\bf R^M}} [\ov(\frac{\ox}{\e}+\beta(z))-\ov(\frac{\ox}{\e})-\la \n_y \ov(\frac{\ox}{\e}), \beta(z) \ra] dq_0(z) + O(\e). 
\end{equation}
(ii)  If $dq_0(z)$ satisfies (\ref{radon}) with $\gamma=1$,  
\begin{equation}\label{est2}
	 \e^{\a} \int_{{\bf R^M}} [\ov(\frac{\ox+\beta(z)}{\e})-\ov(\frac{\ox}{\e})] dq(z)
	=\int_{{\bf R^M}} [\ov(\frac{\ox}{\e}+\beta(z))   -\ov(\frac{\ox}{\e})    ] dq_0 (z)+ O(\e).
\end{equation}
\end{theorem} 

$Proof.$ (i)
From  (\ref{beta}) (i.e. $\e^{-1}\beta(z)=\beta(\frac{z}{\e})$),   we have 
$$
	\e^{\a}\int_{{\bf R^M}} [\ov(\frac{\ox+\beta(z)}{\e})-\ov(\frac{\ox}{\e})-{\bf 1}_{|z|\leq 1}\la \e^{-1}\n_y \ov(\frac{\ox}{\e}), \beta(z) \ra ]dq(z)\qquad\quad
$$
$$
	= \e^{\a} \int_{{\bf R^M}} [\ov(\frac{\ox}{\e}+  \beta(\frac{z}{\e}) )-\ov(\frac{\ox}{\e})-{\bf 1}_{|z|\leq 1}\la \n_y \ov(\frac{\ox}{\e}), \beta(\frac{z}{\e}) \ra ]dq(z)
$$
$$\qquad
	= \int_{{\bf R^M}} [\ov(\frac{\ox}{\e}+  \beta(z') )-\ov(\frac{\ox}{\e})-{\bf 1}_{|\e z'|\leq 1}\la \n_y \ov(\frac{\ox}{\e}), \beta(z') \ra]
	 \e^{M+\a} q(\e z') dz'. 
$$
Then, by the condition (A), (\ref{radon}) with $\gamma=2$, 
$$
	| \int_{{\bf R^M}} [\ov(\frac{\ox}{\e}+  \beta(z') )-\ov(\frac{\ox}{\e})-{\bf 1}_{|\e z'|\leq 1}\la \n_y \ov(\frac{\ox}{\e}), \beta(z') \ra]
	 \e^{M+\a} q(\e z') dz'\qquad\qquad\qquad
$$
$$\qquad\qquad\qquad\qquad
	- \int_{{\bf R^M}} [\ov(\frac{\ox}{\e}+\beta(z))-\ov(\frac{\ox}{\e})- \la \n_y \ov(\frac{\ox}{\e}), \beta(z) \ra] dq_{0}(z)|
$$
$$\qquad
	\leq C  \int_{|z|\leq 1} [\ov(\frac{\ox}{\e}+\beta(z))-\ov(\frac{\ox}{\e})- \la \n_y \ov(\frac{\ox}{\e}), \beta(z) \ra] |\e^{M+\a} q(\e z) - q_{0}(z)| dz 
$$
$$\qquad\qquad\qquad
	+ C  \int_{|z|> 1} |z| |\e^{M+\a} q(\e z) - q_{0}(z)| dz  
$$
$$
	\leq C'  (\int_{|z|\leq 1} |z|^2 |\e^{M+\a} q(\e z) - q_{0}(z)| dz + \int_{|z|> 1} |z| |\e^{M+\a} q(\e z) - q_{0}(z)| dz)=O(\e),
$$
where we used Lemma 4.4 (i) to have the last estimate. 

(ii) The proof is similar to that of (i), while we use (\ref{radon}) with $\gamma=1$ and Lemma 4.4 (ii). \\

Now, we state our first main result of the paper.\\

{\bf Theorem 4.6.$\quad$}
\begin{theorem}  Let $u_{\e}$ be the solution of (\ref{homo})-(\ref{bc}).  Assume that (\ref{beta}), (\ref{radon}) (with $\gamma=2$), (\ref{holder1}), (\ref{holder2}),  and  the conditions (A) and (B) hold. Assume also that $dq_0(z)$ defined in (\ref{q0}) satisfies (\ref{radon}) with $\gamma=2$. 
  Then, there exists a unique function 
$$
	\ou(x)=\lim_{\e\to 0}u_{\e}(x) \quad  \forall x\in {\bf R^N},
$$ which is a unique viscosity solution of 
\begin{equation}\label{effect}
	 \ou(x)+\oI_1[I_1[\ou](x)] = 0 \qquad \hbox{in}\quad \Omega, 
\end{equation}
and (\ref{bc}), where $\oI_1$ is given by (\ref{oI}) with 
$$
	I_1[\ou](x)= \int_{{\bf R^M}} 
	[\ou(x+\beta(z))-\ou(x)- {\bf 1}_{|z|\leq 1}\la \n \ou(x),\beta(z) \ra] dq(z). 
$$
\end{theorem}

$Proof$ $of$ $Theorem$ $ 4.6.$ We use the perturbed test function method introduced in \cite{ev1} (see \cite{lpv}, too), which is now standard to solve homogenization problems in the framework of  viscosity solutions. 
Here, we have to take an extra care to treat the difference between the original L{\'e}vy measure $dq(z)$ and the rescaled measure $dq_0(z)$ in the cell problem (\ref{cell}) (and (\ref{approxcell1})-(\ref{approxcell2})).  Let 
$$
	u^{\ast}(x)=\lim\sup_{\e\to 0,y\to x} u_{\e}(y), \quad u_{\ast}(x)=\lim\inf_{\e\to 0,y\to x} u_{\e}(y) \quad \forall x\in {\bf R^N}.
$$ 
In the following, we devide our argument in two steps.\\
(Step 1.) We  show that $u^{\ast}$ is a subsolution of (\ref{effect}). By assuming that $u^{\ast}$ is not the subsolution of (\ref{effect}), we shall get a  contradiction. So, assume that for a function $\phi(x)\in C^2({\bf R^N})$, $u^{\ast}-\phi$  
takes a global strict maximum at $\ox$, $u^{\ast}(\ox)=\phi(\ox)$, and for some $\gamma>0$, the following holds. 
$$
	\phi(\ox)+ \oI_1[ \int_{{\bf R^M}} [\phi(\ox+\beta(z))-\phi(\ox)- {\bf 1}_{|z|\leq 1}\la \n \phi(\ox),\beta(z) \ra] dq(z)] =3\gamma>0.  
$$
Then, from the continuities of $\oI_1$ and $\phi$, for $r>0$ small enough 
\begin{equation}\label{Br}
	\phi(x)+\oI_1[ I_1[\phi](x)]>2\gamma  \quad \hbox{in}\quad B_r(\ox), 
\end{equation}
where $I_1[\phi](x)=\int_{{\bf R^M}} [\phi(x+\beta(z))-\phi(x)- {\bf 1}_{|z|\leq 1}\la \n \phi(x),\beta(z) \ra] dq(z)$. 

From (\ref{approxcell2}), for $\delta>0$ and $I_1=I_1[\phi](\ox)$, we know that there exists a periodic, semiconcave, Lipschitz continuous function $\ov$ which satisfies
\begin{equation}\label{uv}
	d_{I_1[\phi](\ox)}-a(y)\int_{{\bf R^M}} 
	[\ov(y+\beta(z))-\ov(y)- \la \n \ov(y),\beta(z) \ra] dq_0(z) 
\end{equation}
$$\qquad\qquad\qquad\qquad\qquad\qquad\qquad
	-a(y)I_1[\phi](\ox)-f(y)\geq -\frac{\delta}{2} \qquad \hbox{in}\quad {\bf T^N}.
$$
We claim the following.

{\bf Lemma 4.7.$\quad$}
\begin{theorem}  
Let $\phi_{\e}(x)=\phi(x)+\e^{\a}\ov(\frac{x}{\e})$. The function $\phi_{\e}$ is a viscosity supersolution of 
\begin{equation}\label{Phie}
	\phi_{\e}(x)-a(\frac{x}{\e})\int_{{\bf R^M}} 
	[\phi_{\e}(x+\beta(z))-\phi_{\e}(x)- {\bf 1}_{|z|\leq 1}\la \n \phi_{\e}(x),\beta(z) \ra] dq(z) 
\end{equation}
$$
	\qquad\qquad\qquad\qquad\qquad\qquad
	 -f(\frac{x}{\e})\geq \gamma \quad \hbox{in} \quad B_r(\ox),
$$
where the L{\'e}vy density $dq(z)$ is the one in (\ref{homo}). 
\end{theorem}  
$Proof$ $of$ $Lemma$ $4.7.$ 
To confirm (\ref{Phie}) in the sense of viscosity solutions, assume that for some $\psi\in C^2({\bf R^N})$, $\phi_{\e}-\psi$ takes a strict minimum at 
$x=x'$ and $\phi_{\e}(x')=\psi(x')$. From Definition B in \S 6, we must show 
$$
	\phi_{\e}(x')-a(\frac{x'}{\e})\int_{{\bf R^M}} 
	[\psi(x'+\beta(z))-\psi(x')- {\bf 1}_{|z|\leq 1}\la \n \psi(x'),\beta(z) \ra] dq(z) 
$$
\begin{equation}\label{visco}
	\qquad\qquad\qquad\qquad\qquad\qquad
	 -f(\frac{x'}{\e})\geq \gamma. 
\end{equation}
Since $-(\phi_{\e}-\psi)$ is semiconvex, from Lemma 4.3, we can take a sequence $x'_m\in \Omega$ such that $x'_m\to x'$ as $m\to \infty$, $\phi_{\e}-\psi$ is twice differentiable at $x'_m$, $\n^2(\phi_{\e}-\psi)(x'_m)\geq O$, $|\n(\phi_{\e}-\psi)(x'_m)|\leq \frac{1}{m}$. And by putting $p_m=\n(\phi_{\e}-\psi)(x'_m)$, 
$(\phi_{\e}-\psi)(x)-\la p_m,x \ra$ takes a minimum at $x'_m$. Put $\psi_m(x)=\psi(x)+\la p_m,x \ra$. To see (\ref{visco}), we first prove 
$$
	\phi_{\e}(x'_m)-a(\frac{x'_m}{\e})\int_{{\bf R^M}} 
	[\psi_m(x'_m+\beta(z))-\psi_m(x'_m)- {\bf 1}_{|z|\leq 1}\la \n \psi_m(x'_m),\beta(z) \ra] dq(z) 
$$
\begin{equation}\label{viscom}
	\qquad\qquad\qquad\qquad\qquad\qquad
	 -f(\frac{x'_m}{\e})\geq \gamma, 
\end{equation}
 for any  $m\in {\bf N}$ large enough.  
By remarking that $\phi_{\e}-\psi_m$ is twice differentiable at $x'_m$, that $\psi_m\in C^2$, we know that $\phi_{\e}$ is twice differentiable at $x'_m$, and thus 
 $\phi_{\e}(x_m'+\beta(z))-\phi_{\e}(x_m')- {\bf 1}_{|z|\leq 1}\la \n \phi_{\e}(x_m'),\beta(z) \ra$$\in L^1({\bf R^M},dq(z))$ (we used (\ref{radon})). We can show that 
$$
	\phi_{\e}(x_m')-a(\frac{x_m'}{\e})\int_{{\bf R^M}} 
	[\phi_{\e}(x_m'+\beta(z))-\phi_{\e}(x_m')- {\bf 1}_{|z|\leq 1}\la \n \phi_{\e}(x_m'),\beta(z) \ra] dq(z) 
$$
\begin{equation}\label{Phie'}
	\qquad\qquad\qquad\qquad\qquad\qquad
	 -f(\frac{x_m'}{\e})\geq \gamma, 
\end{equation}
in the classical sense, for any  $m\in {\bf N}$ large enough. 
 To see (\ref{Phie'}),  we use Lemma 4.5 (i) (\ref{est1}) for $\Psi=\phi-\psi_m$, $\ox=x'_m$ to have 
$$
	\e^{\a}\int_{{\bf R^M}} [\ov(\frac{x_m'+\beta(z)}{\e})-\ov(\frac{x_m'}{\e})-{\bf 1}_{|z|\leq 1}\la \e^{-1}\n_y \ov(\frac{x_m'}{\e}), \beta(z) \ra ]dq(z)
$$
$$
	=  \int_{{\bf R^M}} [\ov(\frac{x_m'}{\e}+\beta(z))-\ov(\frac{x_m'}{\e})- \la \n_y \ov(\frac{x_m'}{\e}), \beta(z) \ra] dq_{0}(z) + O(\e). 
$$
Thus, from (\ref{uv}), for $y=\frac{x'_m}{\e}$, $\e>0$ small enough, 
$$
	d_{I_1[\phi](\ox)}-a(\frac{x'_m}{\e})\e^{\a}\int_{{\bf R^M}} [\ov(\frac{x_m'+\beta(z)}{\e})-\ov(\frac{x_m'}{\e})
	-{\bf 1}_{|z|\leq 1}\la \e^{-1}\n_y \ov(\frac{x_m'}{\e}), \beta(z) \ra ]dq(z)
$$
$$\qquad\qquad\qquad\qquad\qquad
	-a(\frac{x'_m}{\e})I_1[\phi](\ox)-f(\frac{x'_m}{\e})\geq -\delta. 
$$

We introduce this into  (\ref{Br}) (for $x=x_m'$$\in B_r(\ox)$):
$$
	\phi(x_m')+\oI_1[\int_{{\bf R^M}}  [\phi(x_m'+\beta(z)) -\phi(x_m') - {\bf 1}_{|z|\leq 1}  \la \n\phi(x_m'), \beta(z) \ra  ]dq(z)   ]>2\gamma.
$$
 By taking $\e>0$, $\delta>0$ small enough so that $\delta+ |\e^{\a}\ov(\frac{x'_m}{\e})|$$\leq \frac{\gamma}{4}$, 
by remarking that $d_{I_1[\phi](\ox)}$ $=-\oI_1(I_1[\phi](\ox))$, from the continuities of $\oI_1$, $\phi$,  for $r>0$ small enough 
we get 
$$
	\phi(x'_m) + \e^{\a}\ov(\frac{x'_m}{\e})-a(\frac{x'_m}{\e})\int_{{\bf R^M}} [(\phi(x'_m+\beta(z))+\e^{\a}\ov(\frac{x'_m+\beta(z)}{\e})  )
$$
$$
	-(\phi(x'_m)+\e^{\a}\ov(\frac{x'_m}{\e}))
	-{\bf 1}_{|z|\leq 1}\la \n \phi(x'_m) + \e^{\a-1}\n_y \ov(\frac{x'_m}{\e}), \beta(z) \ra] dq(z)
$$
$$\qquad\qquad\qquad\qquad\qquad\qquad\qquad\qquad\qquad\qquad\qquad\qquad
	-f(\frac{x'_m}{\e})\geq \gamma. 
$$
Thus, (\ref{Phie'}) is proved. 
From $\n\phi_{\e}(x'_m)$$= \n\psi_m(x'_m)$ and 
$$
	(\phi_{\e}-\psi_m)(x'_m)\leq (\phi_{\e}-\psi_m)(x'_m+ \beta(z)) \quad \forall z\in {\bf R^M}, 
$$
(\ref{Phie'}) leads to (\ref{viscom}): 
$$
	\phi_{\e}(x'_m)-a(\frac{x'_m}{\e})\int_{{\bf R^M}} 
	[\psi_{m}(x'_m+\beta(z))-\psi_{m}(x'_m)- {\bf 1}_{|z|\leq 1}\la \n \psi_{m}(x'_m),\beta(z) \ra] dq(z) 
$$
$$
	\qquad\qquad\qquad\qquad\qquad\qquad
	 -f(\frac{x'_m}{\e})\geq \gamma. 
$$
From (\ref{viscom}), since $|p_m|\leq \frac{1}{m}$, and since 
$$
	\psi_{m}(x'_m+\beta(z))-\psi_{m}(x'_m)- {\bf 1}_{|z|\leq 1}\la \n \psi_{m}(x'_m),\beta(z) \ra
	\qquad\qquad\qquad\qquad\qquad
$$
$$\qquad\qquad
	\to \psi(x'+\beta(z))-\psi(x')- {\bf 1}_{|z|\leq 1}\la \n \psi(x'),\beta(z) \ra\in L^1({\bf R^M},dq(z))
$$
as $m\to \infty$, we have shown (\ref{visco}) in Lemma 4.7.  \\

We continue the proof of Theorem 4.6.  
Now, the comparison principle for  (\ref{homo}) and (\ref{Phie}) leads
$$
	\sup_{x\in U_{r}(\ox)}\{u_{\e}(x)-\phi_{\e}(x)\}\leq \sup_{x\in U_{r}(\ox)^c}\{u_{\e}(x)-\phi_{\e}(x)\}+\gamma. 
$$
By letting $\e\to 0$, since $\gamma>0$ is arbitrary, 
$$
	\sup_{x\in U_{r}(\ox)}\{u^{\ast}(x)-\phi(x)\}\leq \sup_{x\in U_{r}(\ox)^c}\{u^{\ast}(x)-\phi(x)\}. 
$$
However, this contradicts to the fact that $\ox$ is the strict global maximum of $u^{\ast}-\phi$.
 Therefore, $u^{\ast}$ must be a viscosity subsolution of (\ref{effect}). \\
(Step 2.) By the parallel argument, we can prove that $u_{\ast}$ is a viscosity supersolution of (\ref{effect}). Now, from the definition of $u_{\ast}$ and $u^{\ast}$,  we have
$$
	u_{\ast}\leq u_{\e}\leq u^{\ast} \quad \forall \e>0. 
$$
From  
the comparison principle for the viscosity solution of  (\ref{effect})-(\ref{bc}) in Lemma 4.2, we have 
$$
	u^{\ast}\leq u_{\ast} \quad \hbox{in} \quad \overline{\Omega}. 
$$ Thus, there exists a limit: $\ou=\lim_{\e\to 0} $$u_{\e}$$=u_{\ast}=u^{\ast}$ which is the unique viscosity solution of (\ref{effect})-(\ref{bc}). \\

Our second result is the following. \\

{\bf Theorem 4.8.$\quad$}
\begin{theorem}  Let $u_{\e}$ be the solution of (\ref{homo2})-(\ref{bc}).  Assume that (\ref{beta}), (\ref{radon}) (with $\gamma=1$), (\ref{holder1}), (\ref{holder2}) hold,  and  that the conditions (A) and (B) hold. Assume also that $dq_0(z)$ defined in (\ref{q0}) satisfies (\ref{radon}) with $\gamma=1$. 
  Then, there exists a unique function 
$$
	\ou(x)=\lim_{\e\to 0}u_{\e}(x) \quad  \forall x\in {\bf R^N},
$$ which is a unique viscosity solution of 
\begin{equation}\label{effect2}
	 \ou(x)+\oI_2[I_2[\ou](x)] = 0 \qquad \hbox{in}\quad \Omega, 
\end{equation}
and (\ref{bc}), where $\oI_2$ is given by (\ref{oI}) and 
$$
	I_2[\ou](x)= \int_{{\bf R^M}} [\ou(x+\beta(z))-\ou(x)] dq(z). 
$$
\end{theorem} 

$Proof$ $of$ $Theorem$ $ 4.8.$ The proof is similar to that of Theorem 4.6 (in fact, it is simpler because there is no term $ {\bf 1}_{|z|\leq 1}\la \n u(x),\beta(z) \ra$ in the integral). We use Lemma 4.5 (ii) instead of (i). \\

{\bf Corollary 4.9.$\quad$}
\begin{theorem}  (i) Let $u_{\e}$ be the solution of (\ref{homo})-(\ref{bc}).  Assume that (\ref{holder1}), (\ref{holder2}) hold, and that  $dq_0(z)$ and $\beta(z)$ are given either one of the following : Example 1 with $\a\in (1,2)$, 
 Examples 2 and 3 with $\a\in (1,2)$, and Exmple 4 with $\a\in (1,2)$.  Then, there exists a unique function 
$$
	\ou(x)=\lim_{\e\to 0}u_{\e}(x) \quad  \forall x\in {\bf R^N},
$$ which is a unique viscosity solution of (\ref{effect})-(\ref{bc}).
	 \\
(ii) Let $u_{\e}$ be the solution of (\ref{homo2})-(\ref{bc}).  Assume that (\ref{holder1}), (\ref{holder2}) hold, and that  $dq_0(z)$ and $\beta(z)$ are given either one of the following : Example 1 with $\a\in (0,1)$, 
 Example 3 with $\a\in (0,1)$, and Exmple 4 with $\a\in (0,1)$. Then, there exists a unique function 
$$
	\ou(x)=\lim_{\e\to 0}u_{\e}(x) \quad  \forall x\in {\bf R^N},
$$ which is a unique viscosity solution of (\ref{effect2})-(\ref{bc}).
\end{theorem} 
$Proof.$ The claims follows from Corollary 3.3, Theorems 4.6 and 4.8.\\

{\bf Remark 4.2.} 
 The present argument can  be generalized to the following type of the homogenization problem : 
$$
u_{\e}(x)+ 
	\sup_{\tilde{\a}\in \mathcal{A}}\{-a(\frac{x}{\e})\int_{{\bf R^M}}  [u_{\e}(x+\beta(z,\tilde{\a}))-u_{\e}(x)\qquad\qquad\qquad\qquad
$$
$$\qquad\qquad\qquad\qquad- {\bf 1}_{|z|\leq 1}\la \n u_{\e}(x),\beta(z,\tilde{\a}) \ra] dq(z) 
	-f(\frac{x}{\e},\tilde{\a})\}=0\quad \hbox{in}\quad \Omega, 
$$
with (\ref{bc}), where $\mathcal{A}$ is a compact metric set (control set), $\beta(z,\a)$ is a continuous function in ${\bf R^M}\times\mathcal{A}$ with values in 
${\bf R^N}$ satisfying  (\ref{beta}) uniformly in $\mathcal{A}$, $f(y,\a)$ is a real valued continuous function in ${\bf T^N}\times\mathcal{A}$ 
 satisfying  (\ref{holder2}) uniformly in $\mathcal{A}$. We leave the detail to the readers.\\

\section{A nonlinear problem.}

$\quad$In this section, we show how the present method can apply to more general nonlinear 
problems. We consider Example 5 in \S1. Let $u_{\e}$ be the unique viscosity solution of (\ref{ex4}). 

Assume that there exist two 
positive numbers $\a_l\in  (0,2)$ $(l=1,2)$, subsets $S_0^l \subset S^l= \hbox{supp}(dq_l(z))$ $(l=1,2)$, and positive functions $q_{0}^l (z)$  $(l=1,2)$ such that 
the assumption (A) is satisfied: 
\begin{equation}\label{limit2}
	\lim_{\e\to 0} q_l(\e z)\e^{l+\a_l} dz=q_0^l(z) dz\qquad \forall z\in S_0^l; \quad =0 dz\qquad \forall z\in {\bf R}^l / S_0^l
	,\quad l=1,2,
\end{equation}
\begin{equation}\label{bound2}
	|\e^{l+\a_l}q_l(\e z)|\leq C |z|^{-(l+\a_l)}  \quad \forall \e\in (0,1), \quad  \forall z\in {\bf R}^l. 
\end{equation}
where $dq_l(z)=q_l(z)dz$ ($l=1,2$), and $C>0$ is a constant.  
 We define the following new measures : 
\begin{equation}\label{dq0}
	dq_{0}^l (z)=  q_{0}^l (z)dz \quad \forall z\in S_0^l ;\quad =0 dz\quad \forall z\in  {\bf R}^l / S_0^l, \quad l=1,2. 
\end{equation}
 
 Here, we  further assume that $\a_1=\a_2=\a$ (otherwise, a different problem which does not concern with the present interest of the nonlocal problem arises).  
We use the formal asymptotic expansion : 
\begin{equation}\label{nformal}
	u_{\e}(x)=\overline{u}(x)+  \e^{\a} v(\frac{x_1}{\e},\frac{x_2}{\e},\frac{x_3}{\e})  \qquad x\in {\bf R^3}, 
\end{equation} 
and get the following ergodic cell problem. For any given $I'$, $I''\in {\bf R}$, find a unique number $d_{I',I''}$ with which the following problem has a periodic viscosity 
solution $v$: 
$$
	d_{I',I''}+\max\{ -a(y)\int_{{\bf R} } [v(y+\beta_1(z'))-v(y)-\left\langle \beta_1(z'),\n v(y) \right\rangle] dq_0^{1}(z') 
$$
$$
	-a(y)I', 
	 -a(y)\int_{{\bf R^2}} [v(y+\beta_2(z''))-v(y)-\left\langle \beta_2(z''),\n v(y) \right\rangle] dq_0^{2}(z'')
$$
\begin{equation}\label{ncell}
\qquad\qquad\qquad\qquad\qquad\qquad\qquad-a(y)I''
	\}-f(y)=0 
	\quad \hbox{in}\quad {\bf T^3}, 
\end{equation}
where 
$$
	I'=I'[\ou](x)=\int_{{\bf R}} [\ou(x+\beta_1(z'))-\ou(x)-{\bf 1}_{|z'|\leq 1}\left\langle \beta_1(z'),\n \ou(x) \right\rangle] dq_{1}(z'),
$$ 
$$
	I''=I''[\ou](x)=\int_{{\bf R^2}} [\ou(x+\beta_2(z''))-\ou(x)-{\bf 1}_{|z''|\leq 1}\left\langle \beta_2(z''),\n \ou(x) \right\rangle] dq_{2}(z'').
$$ 
As in \S 3, the existence of the unique number $d_{I',I''}$ in (\ref{ncell}) comes from the SMP of the integro-differential equation: 
$$
	H(y,\n v)+ \max\{ -\int_{{\bf R}} [v(y+\beta_1(z'))-v(y)-\left\langle \beta_1(z'),\n v(y) \right\rangle] dq_0^{1}(z'),
$$
\begin{equation}\label{couple}
	 -\int_{{\bf R^2}} [v(y+\beta_2(z''))-v(y)-\left\langle \beta_2(z''),\n v(y) \right\rangle] dq_0^{2}(z'') \}=0 \quad \hbox{in}\quad {\bf T^3}. 
\end{equation}
In order to establish  the SMP for (\ref{couple}), we need to generalize the condition (B) of Theorem 2.1 to the following.  \\

(B') For any two points $y$, $y'\in {\bf T^3}$, there exist a finite number of points $y_1$, ..., $y_m\in {\bf T^N}$ such that $y_1=y$, $y_m=y'$, and for any 
$m$ positive numbers $\e_i>0$ $(1\leq i\leq m)$, we can take subsets 
 $J_i$ ($1\leq \forall i\leq m$) either $J_i\subset S_0^1$ or $J_i\subset S_0^2$, such that  if $J_i\subset S_0^l$ ($l=1,2$),
$$
	 \int_{J_i} 1 dq_0^l (z)> 0; \quad 
	y_i+\b_l(z)\in B_{\e_i}(y_{i+1}) \quad \forall z\in J_i,
$$
for any $1\leq i\leq m$. \\

{\bf Theorem 5.1.$\quad$}
\begin{theorem}  Let $u\in USC({\bf R^3})$  be a viscosity subsolution of (\ref{couple}). 
Assume that $\beta_l$ ($l=1,2$) satisfy (\ref{beta}), that $dq_0^l$  ($l=1,2$) satisfy (\ref{radon}) and  the condition (B'), and that  (\ref{H})  holds. 
 If $u$ attains a maximum at $\oy$ in ${\bf T^3}$, then $u$ is constant in ${\bf T^3}$. 
\end{theorem} 
The proof of Theorem 5.1 is similar to Theorem 2.1, which we do not reproduce here. By using Theorem 5.1, 
the existence of the unique number $d_{I',I''}$ in (\ref{ncell}) can be shown by using a similar argument in \S 3. 
In this way, we can define
$$
	\oI(I',I'')=-d_{I',I''} \quad \forall (I',I'')\in {\bf R^2}. 
$$
Then, the effective integro-differential equation for $\ou$$=\lim_{\e\to 0}u_{\e}$ is the following: 
$$
	\ou+ \oI(I'[\ou](x),I''[\ou](x))=0 \quad x\in \Omega,
$$
associated with the Dirichlet condition (\ref{bc}), where $I'[\ou](x)$ and $I''[\ou](x)$ are given before. This formal argument can be confirmed by the perturbed test function method used in \S 4. Since the argument is similar,   we just show the direction and do not enter in detail here. \\

\section{Appendix.} 
$\quad$ In this section, by following \cite{ardef}, we  note three types of  equivalent definitions of the viscosity solutions for a class of integro-differential equations,  which includes (\ref{homo}). The comparison and the existence of viscosity solutions in this framework are found in  \cite{alvtou}, \cite{arcorrig},  \cite{arcomp}, O. Alvarez and A. Tourin \cite{alvtou}, G. Barles, R. Buckdahn, and E. Pardoux \cite{bbp}, 
  G. Barles and C. Imbert \cite{bi}, and the references there in.  The equivalence of these definitions was shown in \cite{ardef}. We consider the following problem. 
$$
	F(x,u(x),\n u(x),\n^2 u({x}))
	- \int_{\mathbf{R^M}} [u({x}+\beta(z))-u({x})\qquad\qquad\qquad\qquad
$$
\begin{equation}\label{problem}\qquad\qquad\qquad\qquad
	-{\bf 1}_{|z|\leq 1}\la \beta(z),\n u({x})\ra] dq(z) = 0 \quad \hbox{in}\quad \Omega, 
\end{equation}
where $F$ is a real valued continuous function defined in $\Omega\times {\bf R}\times{\bf R^N}\times{\bf S^N}$, which satisfies the degenerate ellipticity (see \cite{users} for the notion).  
 We say that for $u\in USC({\bf R^N})$   (resp. $(LSC({\bf R^N})$),  $(p,X)\in {\bf R^N}\times {\bf S^N}$ is a subdifferential (resp. superdifferential) of $u$ at $x\in \Omega$ if for any small $\mu>0$ there exists $\nu>0$ such that the folowing holds. 
$$
	u(x+z)-u(x)\leq  (resp. \geq)  \quad \la p,z \ra + \frac{1}{2} \la Xz,z \ra + (resp. -) \mu |z|^2 \quad \forall |z|\leq \nu, \quad z\in {\bf R^N}, 
$$
We denote the set of all subdifferentials (resp. superdifferentials)  of $u\in USC(\bf R^N)$ (resp. $LSC(\bf R^N)$) at $x\in \Omega$ by $J^{2,+}_{\Omega}u(x)$ (resp. $J^{2,-}_{\Omega}u(x)$). 
We say that $(p,X)\in {\bf R^N}\times {\bf S^N}$ belongs to  $\overline{J^{2,+}_{\Omega}}u(x)$ (resp. $\overline{J^{2,-}_{\Omega}}u(x)$), if there exist a sequence 
of points $x_n\in \Omega$ and $(p_n,X_n)\in J^{2,+}_{\Omega}u(x_n)$ (resp. $J^{2,-}_{\Omega}u(x_n)$)  such that $\lim_{n\to \infty} x_n=x$, 
 $\lim_{n\to \infty}  (p_n,X_n)=(p,X)$. \\
From (\ref{beta}),  for $u\in USC(\bf R^N)$ (resp. $LSC(\bf R^N)$), if $(p,X)\in J^{2,+}_{\Omega}u(x)$ (resp. $J^{2,-}_{\Omega}u(x)$), we can take 
 a pair of positive numbers $(\nu,\mu)$ such that 
$$
	u(x+\beta(z))-u(x)\leq  (resp.\geq)   \la p,\beta(z) \ra + \frac{1}{2} \la X\beta(z),\beta(z) \ra +(resp.-) \mu |\beta(z)|^2
$$
\begin{equation}\label{subbeta}\qquad\qquad\qquad\qquad\qquad\qquad\qquad\qquad\qquad\qquad
	 \quad \forall |z|\leq \nu, \quad z\in {\bf R^M}, 
\end{equation}

\begin{definition}{\bf Definition A. (\cite{arnewdef})}  Let $u\in USC({\bf R^N})$ (resp. $LSC({\bf R^N})$). We say that $u$ is a viscosity subsolution (resp. supersolution) of (\ref{problem}), if for any $\hx\in \Omega$, any $(p,X)\in J_{\Omega}^{2,+}u(\hx)$ (resp. $\in J_{\Omega}^{2,-}u(\hx)$), and 
any pair of numbers $(\nu,\mu)$ satisfying  (\ref{subbeta}), 
the following holds 
$$
	F(\hat{x},u(\hat{x}),p,X) 
	-\int_{|z|<\nu} \frac{1}{2}\la(X+(resp.-)2\mu I)\beta(z),\beta(z) \ra dq(z)\qquad\qquad\qquad\qquad\qquad\qquad
$$
\begin{equation}\label{def1}
	- \int_{|z|\geq \nu} [u(\hat{x}+\beta(z))-u(\hat{x})
	-{\bf 1}_{|z|\leq 1} \la \beta(z),p\ra] dq(z) \leq(resp. \geq) 0.
\end{equation}
  If $u$ is both a viscosity subsolution and a viscosity supersolution , it is called a viscosity solution.
\end{definition}

\begin{definition}{\bf Definition B.(\cite{bbp},\cite{bi},\cite{jk})}  Let $u\in USC({\bf R^N})$ (resp. $LSC({\bf R^N})$). We say that $u$  is a viscosity subsolution (resp. supersolution) of (\ref{problem}),  if for any $\hx\in \Omega$ and for any $\phi\in C^2({\bf R^N})$ such that $u(\hx)=\phi(\hx)$ and $ u-\phi $ takes a maximum (resp. minimum) at $\hx$, the following holds
$$
	F(\hat{x},u(\hat{x}),\n \phi(\hat{x}),\n^2 \phi(\hat{x}))
	- \int_{{\bf R^M}} [\phi(\hat{x}+\beta(z))
	-\phi(\hat{x})
$$
\begin{equation}\label{defb1}
	 - {\bf 1}_{|z|\leq 1} \la \beta(z),\n \phi(\hat{x})\ra] dq(z) \leq (resp. \geq) 0.
\end{equation}
If $u$ is both a viscosity subsolution and a viscosity supersolution, it is called a viscosity solution.
\end{definition}

\begin{definition}{\bf Definition C. \cite{ardef}}  Let $u\in USC({\bf R^N})$ (resp. $LSC({\bf R^N})$). We say that $u$  is a viscosity subsolution (resp. supersolution) of (\ref{problem}), if for any $\hx\in \Omega$ and  for any $\phi\in C^2({\bf R^N})$ such that $u(\hx)=\phi(\hx)$ and $ u-\phi $ takes a global maximum (resp. minimum) at $\hx$, 
$$
	h(z)=u(\hat{x}+z)-u(\hat{x})-{\bf 1}_{|z|\leq 1}\la \beta(z),\n \phi(\hat{x})\ra \in L^1(\mathbf{R^M}, dq(z)),
$$ 
and 
\begin{equation}\label{defc1}
	F(\hat{x},u(\hat{x}),\n \phi(\hat{x}),\n^2 \phi(\hat{x}))
	- \int_{z\in \mathbf{R^M}} [u(\hat{x}+\beta(z))\qquad\qquad\qquad\qquad
\end{equation}
$$
	\qquad\qquad\qquad\qquad-u(\hat{x})
	-{\bf 1}_{|z|\leq 1}\la \beta(z),\n \phi(\hat{x})\ra] dq(z) \leq(resp.\geq) 0.
$$
If $u$ is both a viscosity subsolution and a viscosity supersolution, it is called a viscosity solution.
\end{definition}

{\bf Theorem 6.1.$\quad$}
\begin{theorem} 
 The definitions A, B, and C are equivalent.
\end{theorem} 

The claim was proved for the case $M=N$ and $\beta(z)=z$ in \cite{ardef}, and for the case $\b$ depending also in $x\in {\bf R^N}$ in \cite{arcomp}. The present case 
is contained in the latter, which is not so different from the former. Thus, we do not reproduce the proof here. \\

We next modify the above definitions  to treat the following 
\begin{equation}\label{problem2}
	F(x,u(x),\n u(x),\n^2 u({x}))
	-\int_{z\in \mathbf{R^M} } [u(x+\beta(z))
	-u(x)] dq(z)=0 \quad \hbox{in}\quad \Omega, 
\end{equation}
which includes (\ref{homo2}), where $dq(z)$ satisfies (\ref{radon}) with $\gamma=1$.  Remark that 
from (\ref{beta}),  for $u\in USC(\bf R^N)$ (resp. $LSC(\bf R^N)$), if $(p,X)\in J^{2,+}_{\Omega}u(x)$ (resp. $J^{2,-}_{\Omega}u(x)$), for any $\mu>0$, we can take 
  $\nu>0$ such that 
\begin{equation}\label{subbeta2}
	u(x+\beta(z))-u(x)\leq  (resp.\geq)   \la p,\beta(z) \ra + (resp.-) \mu |\beta(z)|^2\quad \forall |z|\leq \nu, \quad z\in {\bf R^M}, 
\end{equation}

\begin{definition}{\bf Definition A'.}  Let $u\in USC({\bf R^N})$ (resp. $LSC({\bf R^N})$). We say that $u$ is a viscosity subsolution (resp. supersolution) of (\ref{problem2}), if for any $\hx\in \Omega$, any $(p,X)\in J_{\Omega}^{2,+}u(\hx)$ (resp. $J_{\Omega}^{2,-}v(\hx)$), and 
any pair of positive numbers $(\nu,\mu)$ satisfying  (\ref{subbeta2}), 
the following holds 
$$
	F(\hat{x},u(\hat{x}),p,X) 
	-\int_{|z|<\nu} \la p+ (resp.-) \mu\beta(z),\beta(z) \ra  dq(z)\qquad\qquad
$$
\begin{equation}\label{def12}
	\qquad\qquad\qquad\qquad
	- \int_{|z|\geq \nu} [u(\hat{x}+\beta(z))-u(\hat{x})] dq(z) \leq(resp.\geq)  0.
\end{equation}
 If $u$ is both a viscosity subsolution and a viscosity supersolution , it is called a viscosity solution.
\end{definition}

\begin{definition}{\bf Definition B'.}  Let $u\in USC({\bf R^N})$ (resp. $LSC({\bf R^N})$). We say that $u$ is a viscosity subsolution (resp. supersolution) of (\ref{problem2}),  if for any $\hx\in \Omega$ and for any $\phi\in C^2({\bf R^N})$ such that $u(\hx)=\phi(\hx)$ and $ u-\phi $ takes a maximum (resp. minimum) at $\hx$, and for any $\nu>0$, 
\begin{equation}\label{defb12}
	 F(\hat{x},u(\hat{x}),\n \phi(\hat{x}),\n^2 \phi(\hat{x}))
	- \int_{{\bf R^M}} [\phi(\hat{x}+\beta(z))
	-\phi(\hat{x})] dq(z)\leq (resp.\geq) 0.
\end{equation}
If $u$ is both a viscosity subsolution and a viscosity supersolution, it is called a viscosity solution.
\end{definition}

\begin{definition}{\bf Definition C'.}  Let $u\in USC({\bf R^N})$ (resp. $LSC({\bf R^N})$). We say that $u$ is a viscosity subsolution (resp. supersolution) of (\ref{problem}), if for any $\hx\in \Omega$ and  for any $\phi\in C^2({\bf R^N})$ such that $u(\hx)=\phi(\hx)$ and $ u-\phi $ takes a global maximum (resp. minimum) at $\hx$, 
$$
	h(z)=u(\hat{x}+z)-u(\hat{x}) \quad\in L^1(\mathbf{R^M}, dq(z)),
$$ 
 and 
\begin{equation}\label{defc12}
	F(\hat{x},u(\hat{x}),\n \phi(\hat{x}),\n^2 \phi(\hat{x}))
	- \int_{z\in \mathbf{R^M}} [u(\hat{x}+\beta(z))-u(\hat{x})
	] dq(z) \leq (resp.\geq)0.
\end{equation}
  If $u$ is both a viscosity subsolution and a viscosity supersolution, it is called a viscosity solution.
\end{definition}

{\bf Theorem 6.2.$\quad$}
\begin{theorem} 
 The definitions A', B', and C' are equivalent.
\end{theorem} 

$Proof.$ The proof of Theorem 6.2 can be done in the same way to \cite{ardef}, and we abbreviate it to avoid the redundancy. \\


\begin{thebibliography}{31}
\bibitem{alvtou}
O. Alvarez and A. Tourin, Viscosity solutions of nonlinear integro-differential equations, Ann. Inst. H. Poincar{\'e}, Anal. Non Lin{\'e}aire, 13(1996), pp. 293-317.
\bibitem{arnewdef}
M. Arisawa, A new definition of viscosity solution for a class of second-order degenerate elliptic 
integro-differential equations, IHP Analyse nonlineaire, 23(5) (2006).
\bibitem{arcorrig}
M. Arisawa, Corrigendum for "A new definition of viscosity solution for a class of second-order degenerate elliptic 
integro-differential equations", IHP Analyse nonlineaire, 24(1) (2006).
\bibitem{arloc}
M. Arisawa, A localization of the L{\'e}vy operators arising in mathematical finances, Proceedings in "Stochastic proesses and applications to mathematical finances", World scientics (2007).
\bibitem{ardef}
M. Arisawa, A remark on the definitions of viscosity solutions for the integro-differential equations with L{\'e}vy operators, Journal des Maths Pures et Appliquées (2008). 
\bibitem{arperiodic}
M. Arisawa, Homogenizations for integro-differential equations with L{\'e}vy operators, Communications in  Partial Differential Equations, 34, (2009), no7, pp  617-624. 
\bibitem{arquasi} 
M. Arisawa, Quasi-periodic and almost periodic homogenizations of integro-differential equations with L{\'e}vy operators, in preparation. 
\bibitem{arcomp}
M. Arisawa, Comparison principles for  integro-differential equations with L{\'e}vy operators 
 - the case of spacial depending jumps -, submitted. 
\bibitem{al} 
M. Arisawa and P.-L. Lions, On ergodic stochastic control. Comm. Partial Differential Equations, 23(1998), no.11-12, pp.2187-2217.
\bibitem{bbp}
G. Barles, R. Buckdahn, and E. Pardoux, Backward stochastic differential equations and integral-partial differential equations, Stochastics Stochastics Rep., 
60(1-2)(1997), pp.57-83.
\bibitem{bi}
G. Barles and C. Imbert, Second order elliptic integro-differential equations, Viscosity solutions theory revisited, to appear in IHP Analyse nonlineaire,
 arxiv:math/0702263v1 [math.AP], 2007.
\bibitem{bp1}
G. Barles and B. Perthame, Exit time problems in optimal control and the vanishing viscosity method. SIAM J. Control Optim. 26 (1988), 1133-1148.
\bibitem{blp}
A. Bensoussan, J.L. Lions, and G. Papanicolaou, Asymptotic analysis for
periodic structures. North-Holland, Amsterdam, 1978.
\bibitem{users}
M.G. Crandall, H. Ishii, and P.-L. Lions, User's guide to viscosity
solutions of second order partial differential equations. Bulletin of the AMS,
vol.27, no. 1 (1992).
\bibitem{ev1}
L.C. Evans, The perturbed test function method for viscosity solutions of
nonlinear P.D.E's. Proc. Roy. Soc. Edinburgh, 111A (1989), pp.359-375.
\bibitem{ev2}
L.C. Evans, Periodic homogeneization of certain fully nonlinear partial
differential equations. Proc. Roy. Soc. Edinburgh, 120 A (1992), pp.245-265.
\bibitem{fs}
W.H. Fleming, and H.M. Soner, Controlled Markov processes and Viscosity solutions (1st edition), 
Springer-Verlag 1992.
\bibitem{gm}
M.G. Garroni and J.L. Menaldi, Second order elliptic integro-differential problems, Research Notes in Mathematics 430, Chapman and Hall/Crc 2002. 
\bibitem{jk}
E. Jacobsen and K. Karlsen,  A "maximum principle for semicontinuous functions" applicable to integro-partial differential equations, NoDEA Nonlinear differential equations Appl., 13(2) (2006).
\bibitem{lpv}
P.-L. Lions, G. Papanicolau, and  S.R.S. Varadhan, Homogeneizations of
Hamilton-Jacobi equations. preprint.







\end{thebibliography}
\end{document}